%% file: cuc.tex
\documentclass[11pt]{article}
\usepackage{psfig,amsmath,amssymb,euscript}
\usepackage{lscape}
\renewcommand{\baselinestretch}{1.6}
\makeatletter
\newcommand{\fsize}{\footnotesize}

\input{QYalias}

\begin{document}
\title{\bf Modelling Multivariate Volatilities via Conditionally
Uncorrelated Components\thanks{Partially supported by an EPSRC
research grant and by NSF grant DMS-0355179.}}
\author{
Jianqing Fan$^{1,2}$
\quad \quad Mingjin Wang$^{2,3}$
 \quad \quad Qiwei Yao$^{2,3}$\\[2ex]
$^1$ Benheim Center of Finance and \\
Department of Operations Research and Financial Engineering\\
Princeton University, Princeton, NJ 08544, USA\\[1ex]
$^2$Department of Statistics, London School of Economics, London, WC2A
2AE, UK\\[1ex]
$^3$ Guanghua School of Management, Peking University, Beijing 100871, China}

\date{}

\maketitle

\begin{abstract}
We propose to model multivariate volatility processes based on the
newly defined conditionally uncorrelated components (CUCs). This model
represents a parsimonious representation for matrix-valued processes.
It is flexible in the sense that we may fit each CUC with any
appropriate univariate volatility model. Computationally it splits
one high-dimensional optimization problem into several lower-dimensional
subproblems. Consistency for the estimated CUCs has been established.
A bootstrap test is proposed for testing the existence
of CUCs. The proposed methodology is illustrated  with both simulated and
real data sets.
\end{abstract}


\newpage

\section{Introduction}

One of the most prolific areas of research in the financial
econometrics literature in last two decades is to model
time-varying volatility of financial returns. Many statistical
models, most designed for univariate data, have been proposed for
this purpose. From the practical point of view, there are at least
two incentives to model several financial returns jointly. First,
time-varying correlations among different securities are important
and useful information for portfolio optimization, asset pricing
and risk management. Secondly, modelling for single security may
be improved by incorporating the relevant information in other
securities. The quest for modelling multivariate processes, which
are often represented by conditional covariance matrices, has
motivated the attempts to extending univariate volatility models to
multivariate cases, aiming for practical and/or statistical
effectiveness. We list some of the endeavors below.

Let $\{ \bX_t \}$ be a vector-valued (return) time series with
\[
E(\bX_t | \calF_{t-1} ) = 0, \quad \quad
\var( \bX_t | \calF_{t-1} ) = \bSigma_t \equiv \big( \sigma_{t,ij} \big),
\]
where $\calF_t$ is the $\sigma$-algebra generated by $\{ \bX_t,
\bX_{t-1}, \cdots \}$, and $\bSigma_t$ is an
$\calF_{t-1}$-measurable $d\times d$ semi-positive definite
matrix. One of the most general multivariate GARCH($p,q$) model is
the BEKK representation (Engle and Kroner~1995)
\begin{eqnarray}
\label{a3}
\bSigma_t = \bC + \sum_{i=1}^p \sum_{j=1}^m \bA_{ij} \bX_{t-i} \bX_{t-i}^\tau
\bA_{ij}^\tau
+ \sum_{i=1}^q \sum_{j=1}^m \bB_{ij} \bSigma_{t-i} \bB_{ij}^\tau,
\end{eqnarray}
where $ \bC, \bA_{ij}, \bB_{ij}$ are  $d\times d$ matrices, and
$\bC$ is positive definite (denoted as $\bC>0$).
Although the form of the above model is quite general especially when
$m$ is reasonably large (Proposition~2.2 of Engle and Kroner 1995), it
suffers from the problems of
overparametrization. Similar to multivariate ARMA models, not all
parameters in model (\ref{a3}) are necessarily
identifiable even when $m=1$.
Overparametrization will also lead to a flat likelihood function, making
statistical inference intrinsically difficult and computationally
troublesome. See, for example, Engle and Kroner~(1995), and Jerez, Casals
and Sotoca~(2001).

To overcome the difficulties due to overparametrization, a dynamic
conditional correlation (DCC) model (Engle 2002, Engle and Sheppard~2001)
has been proposed. It is based on the decomposition
\begin{equation} \label{a4}
\bSigma_t = \bD_t \bR_t \bD_t,
\end{equation}
 where $\bD_t = \diag( \sigma_{t,11}^{1/2},
\cdots, \sigma_{t,dd}^{1/2} )$, $\sigma_{t,ii}$ is the conditional
variance of the $i$-th component of $\bX_t$, and $\bR_t \equiv
(\rho_{t,ij})$ is the conditional correlation
matrix. A simple way to facilitate such a model is to model each
$\sigma_{t,ii}$ with a univariate volatility model and to model
conditional correlation using
a rolling exponential smoothing as follows
\[
\rho_{t,ij} =  \sum_{k=1}^{t-1} \la^k \ve_{t-k,i} \ve_{t-k,j}
\Big/ \Big\{ \sum_{k=1}^{t-1} \la^k \ve_{t-k,i}^2
\; \sum_{k=1}^{t-1} \la^k \ve_{t-k,j}^2 \Big\}^{1/2},
\]
where $\ve_{ti}= X_{ti}/\sigma_{t,ii}^{1/2}$. Even with such a
simple specification, estimation typically involves solving a
high-dimensional optimization problem as, for example, the
Gaussian likelihood function cannot be factorized into several
lower-dimensional functions. To overcome the computational
difficulty, Engle~(2002) proposes a two-step estimation procedure
as follows: first fit each $\sigma_{t,ii}$ in (\ref{a4}) with a
univariate GARCH(1,1) model using the observations on the $i$-th
component of $\bX_t$ only, and then  model the conditional
correlation matrix $\bR_t$ by a simple GARCH(1,1) form
\begin{equation}
\label{a5}
{\bR}_t={\bf S}(1-\theta_1-\theta_2)+\theta_1
({\bve}_{t-1}{\bve}_{t-1}^{\prime})+\theta_2 {\bR}_{t-1},
\end{equation}
and $\bve_t$
is a $d\times 1$ vector of the
standardized residuals obtained in the separate GARCH(1,1) fittings for
the $d$ components of $\bX_t$,
and ${\bf S}$ is
the sample correlation matrix of $\bX_t$. Note there are only two unknown parameters
$\theta_1, \theta_2$
in the dynamical correlation model (\ref{a5}), so it can be easily implemented
even for large or very large $d$. However it
may not provide adequate fitting when the components of $\bX_t$ exhibit
different dynamic correlation structures; see an example of three-dimensional
data set in section~4 below. Furthermore in modelling the volatility for
each component, no attempts are made to extract additional
information from other components.

Alexander (2001) proposes an orthogonal GARCH model which fits
each principal component (PC) with a univariate GARCH model
separately, and treats all PCs as {\sl conditionally} uncorrelated
random variables. Since PCs are only unconditionally uncorrelated,
such a misspecification may lead to non-negligible errors in the
fitting; see, for example, Figure~5 and related discussions in
section~4 below.

Other multivariate volatility models include, for
example, vectorized multivariate GARCH models of Bollerslev, Engle
and Wooldridge~(1988), constant conditional
 correlation
multivariate GARCH models of Bollerslev~(1990),
a multivariate stochastic volatility model of Harvey, Ruiz and
Shephard~(1994),
a generalized
orthogonal GARCH models of van der Weide~(2002),
an easy-to-fit ad hoc
approach of Wang and Yao~(2005); see also a survey in Bauwens, Laurent and
Rombouts~(2003) and the references within.

While all the aforementioned models have their own merits,
each of them has one or more of the three drawbacks; (i)
overparametrization, (ii) computational complication, and (iii) too
simple to catch some
important dynamical structures.

In this paper, we propose a new modelling methodology which
mitigates the above three drawbacks. The basic idea is to assume
that $\bX_t$ is a linear combination of a set of {\sl
conditionally uncorrelated components} (CUCs); see section~2.1
below. One fundamental  difference from the orthogonal GARCH model
is that we use CUCs, instead of PCs, which are genuinely
conditionally uncorrelated. The advantages of the new approach
include: (i)~the CUC decomposition leads to a parsimonious
representation for multivariate volatility (matrix-valued)
processes --- there is no model identification problems, (ii)~it has
the flexibility to model each CUC with any appropriate univariate
volatility models, (iii)~computationally it splits a
high-dimensional optimization problem into several
lower-dimensional subproblems, and (iv)~it allows the volatility
model for one CUC to depend on the lagged value of the other
CUCs.

The idea of using CUCs is similar to the so-called the independent
component analysis (Hyv\"arinen, Karhunen and Oja 2001). However
instead of requiring all the component series are independent with
each other, we only impose a weaker condition that the component
series are conditionally uncorrelated; see (\ref{b1}) below. Of
course the existence of CUCs is also not always guaranteed. We
propose a bootstrap test to assess the feasibility of such an
approach. Our empirical experience shows that for a large number
of practical examples, there is no significant evidence to reject
the hypothesis that the CUCs exist.

Literature on applying independent components  analysis to
financial and economic time series includes, for example, Back and
Weigend (1997), Kiviluoto and Oja (1998), M${\breve {\rm
a}}$l${\breve {\rm a}}$roiu, Kiviluoto and Oja (2000), and van der
Weide (2002). Although our basic idea is somehow similar to van
der Weide~(2002), our approach is completely different.

The rest of the paper is organized as follows. Section~2 contains
a detailed description of the proposed new methodology and the
associated theoretical results. Simulation results are reported in
section~3. Illustrations with real data examples are presented in
section~4. Technical proofs are relegated in appendices.

\section{Methodology}

\subsection{Basic setting}

To simplify the matter concerned, we may assume $\var(\bX_t) =
\bI_d$ --- the $d\times d$ identity matrix. In practice, this
amounts to replacing $\bX_t$ by $\bS^{-1/2}\bX_t$, where
$\bS$ is the sample covariance matrix of $\bX_t$.
We assume that each component of $\bX_t$ is a linear
combination of $d$ conditionally uncorrelated components (CUCs)
$Z_{t1}, \cdots, Z_{td}$ which satisfy the conditions $E(Z_{ti}|
\calF_{t-1} )=0$, Var$(Z_{ti}) =1$, and
\begin{equation} \label{b1}
E(Z_{ti}Z_{tj} | \calF_{t-1} ) = 0, \quad \mbox{for all } i\ne j.
\end{equation}
Put $\bZ_t = (Z_{t1}, \cdots, Z_{td})^\tau$.
The above setting implies that
 \begin{equation} \label{b2}
 \bX_t = \bA \bZ_t, \quad \bZ_t = \bA^\tau \bX_t,
 \end{equation}
for a constant matrix $\bA$. Furthermore, $ \var(\bZ_t) = \bA^\tau
\var(\bX_t) \bA = \bA^\tau \bA =\bI_d$. Hence
$\bA$ is a $d\times d$ orthogonal matrix with ${d\over 2}(d-1)$
free elements. Put
\begin{equation} \label{b3}
\var(\bZ_t|\calF_{t-1}) = \diag( \sigma_{t1}^2, \cdots, \sigma_{td}^2),
\end{equation}
i.e. $\sigma_{tj}^2 = \var(Z_{tj} | \calF_{t-1})$. It is easy to see
that once we
have specified $\sigma_{tj}^2$ -- the volatility of
the $j$-th CUC, for $j=1, \cdots, d$,
volatilities for any portfolios can be deduced accordingly. For
example, for any portfolios $\xi_t = \bb^\tau_1 \bX_t$ and $\eta_t
= \bb^\tau_2 \bX_t$ it holds that
\[
\var(\xi_t | \calF_{t-1}) = \sum_{j=1}^d b_{j1}^2
\, \sigma_{tj}^2, \quad \quad \quad
\cov(\xi_t, \eta_t | \calF_{t-1}) = \sum_{j=1}^d b_{j1} b_{j2}
\, \sigma_{tj}^2.
\]
where $(b_{1j}, \cdots, b_{dj}) = \bb_j^\tau \bA$ $(j=1, 2)$.
Hence, the CUC decomposition (\ref{b2}) facilitates a parsimonious
modelling for $d$-dimensional multivariate volatility process via
$d$ univariate volatility models. In this way, we reduce the
number of parameters involved substantially.

\subsection{Estimation of CUCs}

\subsubsection{Estimation procedure}

By (\ref{b2}), $Z_{tj} = \ba_j^\tau \bX_t$, and $ \ba_1, \cdots,
\ba_d$ are $d$ orthogonal vectors. The goal is to estimate the
orthogonal matrix $ \bA =( \ba_1, \cdots, \ba_d) $. Note the
order of $\ba_1, \cdots, \ba_d$ is arbitrary, and cannot be
identified. Furthermore, $\ba_j$ can be replaced by $-\ba_j$.

Condition (\ref{b1}) is equivalent to
\begin{equation} \label{b5}
\max_{B \in \calB_t } \big| E\{ Z_{ti} Z_{tj} I(B) \} \big| = 0
\end{equation}
for any $\pi$-class $\calB_t \subset \calF_{t-1}$ such that the
$\sigma$-algebra generated by $\calB_t$ is equal to $\calF_{t-1}$
(Theorem~7.1.1 of Chow and Teicher, 1997).
In practice, we use some simple $\calB_t$ for the sake of the
tractability. This leads to choosing  an orthogonal matrix $\bA  =
( \ba_1, \cdots, \ba_d )^\tau$ which minimizes
\begin{equation} \label{b6}
\Psi_n(\bA) \equiv
\sum_{1\le i < j \le d} \; \sup_{B \in \calB,\, 1\le k \le k_0 }\;
 {1 \over n-k} \Big|\ba_i^\tau\Big\{ \sum_{t=k+1}^n \bX_t
\bX_t^\tau I( \bX_{t-k} \in B ) \Big\} \ba_j \Big|,
\end{equation}
where $\calB$ is a collection of subsets in $\RR^d$, $k_0 \ge 1$ is
a prescribed integer. We denote by $\wh \bA = ( \wh\ba_1, \cdots,
\wh\ba_d )^\tau$ the resulting estimator.

% Note that when  $\calB$ consists of only two sets, empty set and the
% whole $d$-dimensional space $\RR^d$, $\Psi_n(\bA)$ is basically the same
% as
% $$
% \sum_{1\le i < j \le d} {1 \over n-1} \Big|\ba_i^\tau\Big\{ \sum_{t=2}^n \bX_t
%  \bX_t^\tau  \Big\} \ba_j \Big|.
% $$
% Hence, $\{\wh \ba_j, j=1, \cdots, d\}$ are the principal components.  In
% other words, our model becomes the orthogonal GARCH model in Alexander
% (2001).

Since the order of $\ba_1, \cdots, \ba_d$ is arbitrary, we measure the
estimation error by
\begin{equation} \label{b7}
D(\wh \bA, \; \bA) = 1 - {1\over d}
\sum_{i=1}^d   \max_{1\le j \le d}  | \ba_i^\tau \wh \ba_j| .
\end{equation}
Note that for any orthogonal matrices $\bA$ and $ \bB$, $D(\bA,
\bB)\ge 0$. Furthermore, if the columns of $\bA$ are obtained from
a permutation of the columns of $\bB$ or their reflections, $D(\bA, \bB) = 0$.
In fact $\Psi_n(\bA) = \Psi_n(\bB)$ if and only if $D(\bA, \bB) = 0$.

In practice, we may let $\calB$ consist of balls with an appropriately
selected radius (such that each ball contains sufficiently many data
points) centered on a grid in the sample space of $\bX_t$.
For example, we may use those observations $\bX_t$ as  the centres of balls such
as at least one of the components of $\bX_t$ is the 10th, the 20th, $\cdots$
the 90th sample percentile of the corresponding component observations.

To overcome the difficulties in handling the constraint $ \bA^\tau\bA = \bI_d$
in solving the above optimization problem,
we reparametrize $\bA$ in terms of the decompositions:
\begin{equation} \label{b8}
\bA = \prod_{1\le i < j \le d} \bE_{ij}(\varphi_{ij}),
\end{equation}
where $\bE_{ij}(\varphi_{ij})$ is obtained from the identity matrix
$\bI_d$ with the following replacements: both the $(i,i)$-th and the
$(j,j)$-th elements are replaced by $\cos \varphi_{ij}$, the $(i,j)$-th
and the $(j,i)$-th elements are replaced, respectively, by $\sin
\varphi_{ij}$ and $-\sin \varphi_{ij}$ (Vilenkin 1968, van der Weide 2002).
Obviously $\bE_{ij}(\varphi_{ij})$ is an orthogonal matrix, so is $\bA$
given in (\ref{b8}). Writing $\bA$ in (\ref{b2}) in the form of
(\ref{b8}), the constrained minimization of (\ref{b6}) over
orthogonal $\bA$ is transformed to an unconstrained minimization
problem over a ${d(d-1)\over 2}\times 1$ vector $\bvarphi =
(\varphi_{12}, \varphi_{13},\cdots, \varphi_{1d}, \varphi_{23}, \cdots,
\varphi_{d-1,d} )^\tau$. This minimization problem is typically
solved by iterative algorithms.
We stop the iteration when $D(\bA_k, \bA_{k+1})$ is
smaller than a prescribed small number, where $\bA_k$ denotes the
value of $\bA$ in the $k$-th iteration, and $D$ is defined as in
(\ref{b7}).

\noindent {\bf Remark 1}. In practice,  we may replace (\ref{b6})
by a weighted version
$$
\Psi_n(\bA) = \sum_{1\le i < j \le d} \; \sup_{B \in \calB,\, 1\le
k \le k_0 }\;
 {1 \over n-k} \Big|\ba_i^\tau\Big\{ { \sum_{t=k+1}^n \bX_t
\bX_t^\tau [I( \bX_{t-k} \in B ) + \varepsilon_0] \over
\sum_{i=k+1}^n [I( \bX_{t-k} \in B )+\varepsilon_0] } \Big\} \ba_j
\Big|,
$$
where $\varepsilon_0$ is a small constant guarding against zero
denominator.  This puts more emphasis on small sets $B$.
Furthermore, the superemum over $k$ in (\ref{b6}) may be replaced
the summation over~$k$.

\subsubsection{Asymptotic properties}

We first introduce two concepts:  mixing which measures the decaying speed of
the auto-dependence for a time series over an increasing time span, and
the Vapnik-$\breve{\mbox{C}}$ervonenkis (or VC) index which measures
the complexity of a collection of sets.

Let $\calF_{i}^j$ be the $\sigma$-algebra generated by $\{\bX_t, i
\leq t \leq j \}$. The  $\beta$-mixing coefficients is defined  as
$$
 \beta(n) = E \left \{ \sup_{ B \in
  \calF_n^\infty} | P(B) - P(B|\calF_{-\infty}^0 ) | \right \}.
$$
(See \S 2.6.1 of Fan and Yao, 2003.)

For an arbitrary set of $n$ points $\{x_1, \cdots, x_n \}$, there are
$2^n$ possible subsets.  Say that $\calB$ picks out a certain subset
from $\{x_1, \cdots, x_n\}$ if this can be formed as a set of the
form $B \cap \{x_1, \cdots, x_n\}$ for a set $B$ in $\calB$. The
collection $\calB$ shatters $\{x_1, \cdots, x_n\}$ if each of its
$2^n$ subsets can be picked out by $\calB$.  The VC-index of $\calB$
refers
to the smallest $n$ for which no set of size $n$ is shattered by
$\calB$. A collection of sets $\calB$ is called a VC-class if its
VC-index is finite.  The collections of sets of rectangles, balls and
their unions are VC-classes. See Chapter 2.6 of van der Vaart and
Wellner (1996) for further discussion on VC-classes.

Under the regularity conditions listed below, the estimator $\wh \bA$
is consistent; see Theorem~1. Its proof is relegated in Appendix A.
\begin{quote}
 (A1) The collection $\calB$ of sets  in $\RR^d$  is a VC-class.

 (A2) The process $\{ \bX_t \}$ is strictly stationary with $E||
\bX_t ||^2 < \infty$, where $||\cdot||$ denotes the Euclidean
norm. Furthermore, the $\beta$-mixing coefficients  $\{\bX_t \}$
satisfy $\beta(n) = O(n^{-b})$ for some $b > 0$.

(A3) There exists a $d\times d$ orthogonal matrix $\bA_0$ which
minimises
\[
\Psi(\bA) \equiv \sum_{1\le i < j \le d}  \sup_{1 \leq k \leq k_0,
B \in \calB} \big| E \{ \ba_i^\tau \bX_t \bX_t^\tau \ba_j
I(\bX_{t-k} \in B) \} \big|.
\]
Furthermore the minimum value of $\Psi$ is obtained at an orthogonal
matrix $\bA$ if and only if $D(\bA, \bA_0) = 0$.

(A4).  $E \| \bX_t \|^{2p} < \infty$ for some $p >2$ and the
$\beta$-mixing coefficient in (A2) holds for $b > p/(p-2)$.

(A5) $\Psi(\bA_0) - \Psi(\bA) \le - a D(\bA, \bA_0)$
for any orthogonal matrix $\bA$ such that $D(\bA, \bA_0)$ is smaller than a
small but fixed constant, where $a > 0$ is a constant.
\end{quote}

\noindent{\bf Remark 2}.  Let $\calH$ be the set consisting of all
$d\times d$ orthogonal matrices.
Then $\calH$ may be partitioned into the equivalent classes defined
by the distance $D$ in (\ref{b7}) as follows: the $D$-distance
between  any two elements within an equivalent class is 0, and the
$D$-distance between
any two elements from different classes is greater than 0.
Let $\calH_D$ be the quotient space $\calH/D$ consisting of those
equivalent classes in $\calH$, i.e.  we treat $\bA$ and $\bB$ as the
same element in $\calH_D$ if and only if $D(\bA, \bB) =0$.
Condition (A3) ensures $\bA_0$ is the unique minimiser
of $\Psi(\bA)$ on $\calH_D$.
In fact both $\Psi(\cdot)$ and $\Psi_n(\cdot)$ are
 Lipschitz continuous on $\calH_D$ with $D$-distance; see Lemma~1 in Appendix~A
below.

\askip

\noindent {\bf Theorem 1}. Let $k_0\ge 1$ be a fixed integer.
Under conditions (A1)--(A3), $D(\wh \bA, \bA_0) \to 0$ almost
surely as $n \to \infty$.  If, in addition, condition (A4) holds,
then
$$
\Psi_n (\bA) - \Psi(\bA) = O_P(n^{-1/2}), \quad \mbox{for any
orthogonal $\bA$.}
$$
Furthermore, $n^{1/2} D(\wh\bA, \bA_0) = O_P(1)$ provided that, in
addition, condition (A5) also holds.

When the CUCs exist, namely $\Psi(\bA_0) = 0$, $\bA_0$ corresponds to the 
transform for the CUCs.  When the CUC does not exist, Theorem 1
continues to hold.  In this case, $\Psi(\bA_0) \not = 0$ and indeed
$\bA_0$ can depend on the $\pi$-class $\cal B$.
In practice, we really do not know whether this condition holds
or not. In that case, our aim becomes naturally to find an
orthogonal transform such that the resulting components are as
less conditionally correlated as possible.  Observe that the
conditional correlation criterion
$$
  \Psi(\bA) = \sum_{1\le i < j \le d}  \sup_{1 \leq k \leq k_0,
   B \in \calB} \big| \mbox{Corr} (\ba_i^\tau \bX_t, \ba_j^T \bX_t |
   \bX_{t-k} \in B ) \big | P( \bX_{t-k} \in B).
$$
Thus, a reasonable criterion is to find an orthogonal transform
$\bA$ to minimize $\Psi(\bA)$.  The following theorem shows that
our estimation method possesses some degrees of robustness and is
better than the principal component transform in terms of
minimizing the conditional correlation criterion $\Psi(\bA)$.

\noindent {\bf Theorem 2}. Let $k_0\ge 1$ be a fixed integer.
Under conditions (A1), (A2),  for any other orthogonal transform
$\hat{\bB}$, we have
$$
   \liminf  \{\Psi(\hat{\bA}) - \Psi(\hat{\bB})\} \leq 0.
$$
%If, in addition, $\bA_0$ is the unique minimizer of $\Psi(\bA)$ on
%the quotient space $\calH_D$, then $D(\hat{\bA}, \bA_0) \to 0$
%almost surely.

Theorem 2 shows for any other orthogonal transform $\hat{\bB}$,
asymptotically, the transformed components have higher conditional
correlation, in terms of $\Psi(\cdot)$, than those transformed by
$\hat{\bA}$.

\subsection{Modelling volatilities for CUCs}

Once the CUCs have been identified, we may fit each
$\sigma_{tj}^2$ with any appropriate univariate volatility model,
for example, a GARCH model, a stochastic volatility model, or any
nonparametric and semiparametric volatility models. As a simple
illustration, we establish below an extended GARCH(1,1) model for
each of $\sigma_{ti}^2$ given in (\ref{b3}).

\subsubsection{Extended GARCH(1,1) models}

We assume, for the $j$-th CUC, $j=1, \cdots, d$,
\begin{equation} \label{b9}
Z_{tj} = \sigma_{tj} \ve_{tj}, \quad \quad \sigma_{tj}^2 = \ga_j +
\sum_{i=1}^d \alpha_{ji} Z_{t-1, i}^2 + \beta_j \sigma_{t-1,j}^2,
\end{equation}
where $ \{\ve_{tj}, \; -\infty < t < \infty\} $ is a sequence of i.i.d.
random variables with mean 0 and
variance 1, $\ve_{tj}$ is independent of $\calF_{t-1}$, $\ga_j
>0$ and $\alpha_j, \alpha_{ji}, \beta_j \ge 0$.
This model contains extra $d-1$ terms $\sum_{i \not = j}
\alpha_{ji} Z_{t-1, i}^2$ from the standard GARCH(1,1) model,
which incorporates the possible association between the $j$-th CUC
and the other CUCs, while the conditional zero-correlation
condition (\ref{b1}) still holds. Such a dependence is termed as
that the $i$-th component (if $\alpha_{ji} \not = 0$) is causal in
variance to the $j$-th component (Engle, Ito and Lin~1991).

In practice, we expect that $\sigma_{tj}^2$ may depend on
$Z_{t-1, i}^2$ only for a small number of $i$'s, including $i=j$, i.e. many
coefficients $\alpha_{ji}$ (for $i\ne j$) may be 0.
Section~2.3.3 below outlines a data-analytic approach for
building such a component-dependent model.

When $\beta_j \in [0, 1)$, (\ref{b9}) implies
\begin{equation} \label{b10}
\sigma_{tj}^2 = \var(Z_{tj} |\calF_{t-1}) = {\ga_j \over 1
-\beta_j} + \sum_{i=1}^d \alpha_{ji} \sum_{k=1}^\infty
\beta_j^{k-1} Z_{t-k,\, i}^2.
\end{equation}
Put $\bZ_t = (Z_{t1}, \cdots, Z_{td})^\tau$. Theorem~2 below gives
a sufficient condition of the existence of stationary solution to
model~(\ref{b9}).

\askip

\noindent {\bf Theorem 3}. (i) The extended GARCH(1,1) model
(\ref{b9}) defines a unique $d$-dimensional strictly stationary
process $\{ \bZ_t \}$ with $E || \bZ_{t}||^2 < \infty$ under the
condition
\begin{equation} \label{b11}
r\cdot \max\{\alpha_{j1}, \cdots, \alpha_{jd} \} + \beta_j < 1,
\quad \quad 1\le j \le d,
\end{equation}
where $r = \max_{1\le j \le d} d_j$, and $d_j$ is the number of non-vanishing
coefficients among $ \alpha_{j1}, \cdots, \alpha_{jd} $.

(ii) Under condition (\ref{b11}), $E(Z_{tj}^2) = 1 $ for all $1\le
j \le d$ if and only if
\begin{equation}
\ga_j = 1 - \beta_j - \sum_{i=1}^d \alpha_{ji} , \quad \quad \quad
1\le j \le d.  \label{b12}
\end{equation}

\askip

The proof of the above theorem is in Appendix B. When
$\alpha_{ji}= 0$ for all $i \not = j$, i.e. each $Z_{tj}$ follows
a standard GARCH(1,1) model, (\ref{b11}) reduces to $\alpha_{jj} +
\beta_j < 1$, which is the necessary and sufficient condition for
the existence of unique strictly stationary solution with finite
second moments for the corresponding GARCH(1,1) model; see Chen
and An (1998).  In practice condition (\ref{b11}) may often be
violated, indicating the likely inappropriateness of GARCH
specification for $\sigma_{tj}^2$. However if we view the right
hand side of (\ref{b10}) as an approximation for $\sigma_{tj}^2$,
such an approximation process is strictly stationary under a weaker
condition $\beta_j <1$. For further discussion of the
approximation point of view, we refer to Penzer, Wang and Yao~(2004).

\subsubsection{quasi-MLE}

To facilitate a likelihood, let us assume hypothetically
that $\ve_{tj}$ in (\ref{b9}) has a density $f(\cdot)$,
which can be the standard normal distribution, generalized Gaussian
distribution and $t$-distribution.  The implied (negative)
log-likelihood function for $\btheta_j \equiv (\alpha_{j1},
\cdots, \alpha_{jd}, \beta_j)^\tau$ is
\begin{equation} \label{b13}
l_j(\btheta_j ) = \sum_{t=\nu+1}^n \big\{
\log \sigma_{tj}(\btheta_j)  - \log f(Z_{tj}/\sigma_{tj}(\btheta_j)) \big\},
\end{equation}
for a given integer $\nu \ge 1$, where $\sigma_{tj}(\btheta_j)^2 =
\var(Z_j | \calF_{t-1})$ is given by $(\ref{b9})$.  By (\ref{b10})
and  (\ref{b12}),
\begin{eqnarray}
\sigma_{tj}(\btheta_j)^2 &=&
 \frac{\gamma}{1 - \beta_j} + \sum_{i=1}^d \alpha_{ji} \sum_{k=1}^\infty
    \beta_j^{k-1} Z_{t-k,i}^2 \nonumber \\
&=& 1 - \frac{ \sum_{i=1}^d \alpha_{ji}}{ 1 -\beta_j}  +
\sum_{i=1}^d \alpha_{ji} \sum_{k=1}^\infty \beta_j^{k-1}
Z_{t-k,i}^2 . \label{b14}
\end{eqnarray}
This form of $\sigma_{tj}(\btheta_j)^2  $ ensures
$\var(Z_{tj})=1$; see Theorem~2(ii). The quasi-maximum likelihood
estimator $\wt \btheta_j$ minimizes (\ref{b13}). In practice, we
let $Z_{ti} \equiv 0$ for all $t\le 0$ on the right hand side of
(\ref{b14}).

\subsubsection{Selection of casual components}

To obtain a parsimonious representation for $\sigma_{tj}^2$, we
may select only those significant $Z_{t-1,i}$ on the RHS of the
second equation in (\ref{b9}). This is particularly important when
the number of components $d$ is large. It may be achieved by using
the ideas for variable selection in regression analysis. Below we
outline such an algorithm based on a combination of the stepwise
addition method and the BIC criterion.

We start with the standard GARCH(1,1) model (i.e. $\alpha_{jj}\ne 0$
and $\alpha_{ji} = 0$ for $j \not = i$). We then add one more $Z_{t-1,i}$
each time which maximizes the (quasi-)likelihood.
More precisely, suppose the model contains
$(k-1)$ terms $Z_{t-1, j_1}, \cdots, Z_{t-1, j_{k-1}}$ already.
 We choose an additional term $Z_{t-1, \ell}$ among
$\ell\not\in \{j, j_1, \cdots, j_{k-1}\}$ which maximizes the
quasi-likelihood function. Note that this is a two-step
maximization problem:  For each given $\ell\not\in \{j, j_1,
\cdots, j_{k-1}\}$, we compute the qMLE $\wt \btheta_j^{(k)}$ for
$\btheta_j^{(k)} \equiv  (\alpha_{jj}, \alpha_{jj_1}, \cdots,
\alpha_{j\ell}, \beta_j)^\tau$ with the constraints
$\alpha_{ji} = 0$, for $i \not \in \{j, j_1, \cdots, j_{k-1}, \ell\}$. We then choose
an $\ell \not \in \{j, j_1, \cdots, j_{k-1} \}$ to minimize
$l_j(\wt \btheta_j^{(k)})$, and denote by $l_j(k)$ the minimum
value and the index of the selected variable $j_k$. Put
\[
{\rm BIC}_j(k) = l_j(k) +  (k+2) \log(n-\nu).
\]
We choose $r_j$ which minimizes BIC$_j(k) $ over $0 \le k \le d$.
Note that $k=0$ corresponds the standard GARCH(1,1) fitting for
$Z_{tj}$.

\subsubsection{LADE}

If CUCs $Z_{tj}$ are known (i.e. $\ba_j$ are known), the asymptotic properties
of qMLE may be derived in the similar manner as Hall and
Yao~(2003). See also Mikosch and Straumann~(2004). For example, the
estimator $\wt \btheta_j$ would suffer from
complicated asymptotic distributions and slow convergence rates if $\ve_{tj}$
is heavy-tailed in the sense that $E(|\ve_{tj}|^4) = \infty$.
On the other hand, a least absolute deviation estimator
based on a log-transformation
is always asymptotically normal with the standard root-$n$
convergence rate; see Peng and Yao (2003).

To construct the LADE with the constraint $\var(Z_{tj})=1$, we
write $\ve_{tj} = v_0 e_{tj}$ in the first equation in (\ref{b9}),
where the median of $e_{tj}^2$ is equal to 1 and $v_0 =
1/{\mbox{STD}}(e_{tj})$. With $\sigma_{tj}(\btheta_j)^2$ expressed
in (\ref{b14}), parameters $\btheta_j$ and $v_0$ are (jointly)
identifiable. Now
\[
\log Z_{tj}^2 - \log \{ \sigma_{tj}(\btheta_j)^2\} - \log v_0^2
= \log (e_{tj}^2).
\]
Since the median of $\log (e_{tj}^2) $ is 0, the true values of the
parameters minimise
\[
E \big|\log Z_{tj}^2 - \log \{ \sigma_{tj}(\btheta_j)^2\} - \log v_0^2\big|.
\]
Therefore we may estimate the
parameters by minimizing
\begin{equation} \label{b15}
\sum_{t=\nu+1}^n
|\log Z_{tj}^2 - \log \{ \sigma_{tj}(\btheta_j)^2\} - \log v_0^2\big|,
\end{equation}
where $\sigma_{tj}(\btheta_j)^2$ is given in (\ref{b14}), with the
part of $a_{ji} = 0$ for the non-casual component in the variance.
So far $\btheta_j$ and $v_0$ are treated as free parameters. The estimators
obtained are root-$n$ consistent.

To make an explicit use of the condition that $\var(\ve_{tj})=1$,
we may estimate parameters $\btheta_j$ as follows.
With the initial estimate $\hat{\btheta}_j^{(0)}$, let $\hat{v}_0$ be the
reciprocal of the sample standard deviation of the residuals $\{
\wt\ve_{tj} \}$, where $\wt\ve_{tj}
=Z_{tj}/\{\sigma_{tj}(\btheta_j^{(0)}) \}$.
With the given $\hat{v}_0$ and $\hat{\btheta}_j^{(0)}$, we can
minimize
$$
\sum_{t=\nu+1}^n  w_t \bigl ( \log Z_{tj}^2 - \log \{
\sigma_{tj}(\btheta_j)^2\} - \log \hat{v}_0^2\bigr )^2,
$$
where $w_t =  |\log Z_{tj}^2 - \log \{
\sigma_{tj}(\hat{\btheta}_j^{(0)})^2\} - \log
\hat{v}_0^2\big|^{-1}$. We may update $\hat{v}_0$ and
iterate further until the estimated
$\btheta_j$ converges. Note that we have used a weighted $L_2$ loss
function to approximate the $L_1$ loss to expedite the computation.

\subsection{Inference based on bootstrapping }

A natural question for the proposed approach is if the CUCs
$Z_{t1}, \cdots, Z_{td}$ exist, although the minimiser $\{ \wh
\ba_j\}$ of (\ref{b6}) always exists. To address this issue
statistically, we may construct a test for the null hypothesis
\[
H_0: \; \bX_t = \bA \bZ_t \quad \mbox{and} \quad
\bZ_t = \diag(\sigma_{t1}, \cdots, \sigma_{td}) \bve_t,
\]
where $\bA^\tau\bA = \bI_d$, $\bve_t = (\ve_{t1}, \cdots,
\ve_{td})^\tau$, $\{ \ve_{t1}\}, \cdots, \{ \ve_{td}\}$ are $d$
independent series, and each of them is a sequence of i.i.d. r.v.s
with mean 0 and variance 1. Note that the null hypothesis above is
a sufficient but not necessary condition for the existence of
CUCs. The independence condition is required to construct a
bootstrap test for this null hypothesis.

Note when $Z_{ti}$ and $Z_{tj}$ are not conditionally
uncorrelated, the left hand side of (\ref{b5}) is equal to
positive constant instead of 0. Therefore, the {\sl large} values
of $\Psi_n(\wh \bA)$ will indicate that the CUCs do  not exist. We
adopt a bootstrap method below to assess how large is large enough
to reject~$H_0$.

If the null hypothesis $H_0$ could not be rejected, we may also
construct confidence sets for the coefficients $\ba_j$ (i.e. the
columns of $\bA$) of the CUCs, and the parameters $\btheta_j$
based on the same bootstrap scheme. Formally confidence sets for
$\btheta_j$ could be constructed based on asymptotic distributions
of, for example, the LADE $\wh \btheta_j$, which may be derived in
the similar manner of Peng and Yao~(2003). However such an
approach is based on the assumption that the CUCs are known (i.e.
the vectors $\ba_j$ are known), and, therefore, fails to take
into account of the errors due to the estimation for $\ba_j$.

Let $\wh \bA =(\wh \ba_1, \cdots, \wh\ba_d)$ be the estimator
derived from minimizing (\ref{b6}). Let $Z_{tj} = \wh \ba_j^\tau
\bX_t$. Let $\wh \btheta_j$ be an estimator
for $\btheta_j$, such as the LADE defined in section~2.3.4.

The bootstrap sampling scheme consists of the three steps below.
\begin{quote}
(i) For $j=1, \cdots, d$, draw $\ve_{tj}^*$, for $-\infty< t \le n$,
by sampling randomly with replacement from  the standardized residuals
$\{\wh \ve_{\nu+1, j}, \cdots , \wh \ve_{nj}\}$ which are obtained
from standardizing the raw residuals
\[
Z_{tj}/\sigma_{tj}(\wh \btheta_j), \quad \quad t=\nu+1, \cdots, n.
\]

(ii) For $j=1, \cdots, d$, draw $Z_{tj}^* = \sigma_{tj}^* \ve_{tj}^*$,
for $-\infty< t \le n$, where
\[
( \sigma_{tj}^*)^2 =1 -  \wh \beta_j - \sum_{i=1}^d \wh
\alpha_{ji}   + \sum_{i=1}^d \wh \alpha_{ji}(Z_{t-1, i}^*)^2 + \wh
\beta_j (\sigma_{t-1,j}^*)^2.
\]

(iii) Let $\bX_t^* = \wh \bA (Z_{t1}^*, \cdots, Z_{td}^*)^\tau$ for $t=1, \cdots,
n$.
\end{quote}

\askip

\noindent {\sl A test for the existence of the CUCs}: Let
$\Psi_n^*(\bA)$ be defined as in (\ref{b6}) with $\{ \bX_t \}$
replaced by $\{ \bX_t^* \}$, and the bootstrap estimator $\bA^*=
(\ba_1^*, \cdots,  \ba_d^*)$ be computed in the same manner as
$\wh \bA$ with $\{ \bX_t \}$ replaced by $\{ \bX_t^* \}$. Note
that the bootstrap sample $\{ \bX_t^* \}$
 is drawn from the model with $\wh \ba_j^\tau \bX_t$ as its {\sl genuine}
CUCs. Hence the conditional
distribution of $\Psi_n^*( \bA^*)$ (given the original sample $\{ \bX_t
\}$) may be taken as an approximation for the distribution of $\Psi_n(\wh
\bA)$ under $H_0$.  Thus we reject
$H_0$ if $\Psi_n(\wh \bA)$ is greater than the $[B\alpha]$-th largest
value of $\Psi_n^*( \bA^*)$ in a replication of the above bootstrap
resampling for $B$ times, where $\alpha \in (0, 1)$ is the size of the
test and $B$ is a large integer.

\askip

\noindent {\sl Confidence sets for ${\bf A}$}: A bootstrap
approximation for an $(1-\alpha)$ confidence set of the
transformation matrix ${\bf A}$ can be constructed  as
\begin{equation}\label{b16}
\{ {\bf A} \, \big| \, D(\hat{\bf A}; {\bf A}) \le c_\alpha , {\bf A}^{\tau}{\bf A}={\bf I}_d \},
\end{equation}
where  $c_\alpha $ is the $[B\alpha]$-th largest value of
$D(\hat{\bf A}; {\bf A}^*)$ in a replication of bootstrap
resampling for $B$ times. Note that when $\bA$ is in the
confidence set, so is $\bB$ if the columns of $\bB$ form a
permutation of the (reflected) columns of $\bA$; see (\ref{b7}).

\askip

\noindent
{\sl Interval estimators for the components of $\wh \btheta_j$}:
A bootstrap confidence interval for any component, say, $\beta_j$
of $\btheta_j$ may be obtained as follows. Repeat the above
bootstrap sampling $B$ times for some large integer $B$, resulting
in bootstrap estimates $ \beta^*_{j1}, \cdots,  \beta^*_{jB}$. An
approximate $(1-\alpha)$ confidence interval for $\beta_j$ is $(
\beta_{j(b_1)}^*, \; \beta_{j(b_2)}^*)$, where $ \beta_{j(i)}^*$
denotes the $i$-th smallest value among $ \beta^*_{j1}, \cdots,
\beta^*_{jB}$, and $ b_1=[B\alpha/2]$ and $b_2=[B(1-\alpha/2)]$.

\section{Simulation}

We conduct a Monte Carlo experiment to illustrate the proposed
CUC-approach. In particular we check the accuracy of the
estimation for the transformation matrix $\bf A$ in (\ref{b2}).

We consider a
CUC-GARCH(1,1) model with $d=3$
\begin{equation} \label{ex1}
{\bf X}_t  =  {\bf A}{\bf Z}_t, \quad \quad
            {\bf Z}_{t}| \calF_{t-1}\; \sim \;  N(0,\;  \diag\{\sigma_{t,1}^2,
\sigma_{t,2}^2, \sigma_{t,3}^2\}),
\end{equation}
where $ \sigma_{t,i}^2  =  \gamma_i+\alpha_i
Z_{t-1,i}^2+\beta_i \sigma_{t-1,i}^2 $, and
\begin{center}
\begin{tabular}{ccc|cccc}
     & {\bf A}  &        &   $i$       &   $\gamma_i$   &  $\alpha_i$
  &  $\beta_i$  \\[0.5ex]\hline
  0  &   0.500  & 0.866  &     1             & 0.02         &   0.08
  &      0.90 \\
  0  &   0.866  & -0.500 &     2             & 0.10         &   0.10
  & 0.80 \\
  -1 &   0      & 0      &     3             & 0.28         &   0.12
  & 0.60 \\
\end{tabular}
\end{center}
It is easy to see that ${\bf A}^\tau{\bf A}={\bf I}_3$ and
$\gamma_i = 1 - \alpha_i - \beta_i$ so that the variances of the
CUCs are 1 [see (\ref{b12})]. Since $\alpha_1 + \beta_1 = 0.98$
is very close to 1, the volatility for the first CUC is highly
persistence. On the contrary, the volatility persistence in the
third component is less  pronounced as $\alpha_3+\beta_3=0.72$
only.

For each of 200 samples with size $n=500$ and 1000 respectively
from the above model, we estimated the transformation matrix $\bA$
by minimizing $\Psi_n({\bf A})$ defined in (\ref{b6}), which was
solved using the proprietary optimization routines in MATLAB. Note
that as far as the estimation of $\bA$ is concerned, two
orthogonal matrices are treated as identical if the $D$-distance
between them is 0; see (\ref{b7}). The coefficients $\alpha_i,
\beta_i$ and $\gamma_i$ were estimated using quasi-MLE based on a
 Gaussian likelihood. The
resulting estimates were summarized in Table~1 and Figure~1.

\begin{table}[htb]
\begin{center}
\caption[Table 1]{Simulation Results: summary statistics of the
errors in estimation}
\begin{tabular}{cc | c c c c c c c}\hline
&& $D(\hat{\bf A}, {\bf A})$ &    $\hat{\alpha}_1$  & $ \hat{\beta}_1$  &  $ \hat{\alpha}_2$  &  $ \hat{\beta}_2$  &  $\hat{\alpha}_3$     &  $\hat{\beta}_3$  \\[0.5ex]\hline
&  mean       &   0.0753      &    0.0719            &      0.8701       &
     0.0865         &    0.7506          &     0.0997            &
0.6189      \\
&  median     &   0.0474      &    0.0705            &      0.8870       &
     0.0830         &    0.7801          &     0.0861            &
0.6445      \\
$n=500$ & STD   &   0.0714      &    0.0300            &      0.0830       &
     0.0469         &    0.1469          &     0.0600            &
0.2017      \\
&  bias       &      -        &   -0.0081            &     -0.0299       &
    -0.0135         &   -0.0494          &    -0.0203            &
0.0189      \\
&  RMSE       &      -        &    0.0303            &      0.0888       &
     0.0484         &    0.1546          &     0.0629            &
0.2022      \\ \hline
&  mean       &   0.0679      &    0.0722            &      0.8921       &
     0.0846         &    0.7751          &     0.0937            &
0.6307      \\
&  median     &   0.0434      &    0.0731            &      0.8999       &
     0.0833         &    0.7956          &     0.0938            &
0.6517      \\
$n=1000$&  STD   &   0.0648      &    0.0224            &      0.0400       &
     0.0346         &    0.1065          &     0.0412            &
0.1634      \\
 & bias       &      -        &   -0.0078            &     -0.0079       &
    -0.0154         &   -0.0249          &    -0.0263            &
0.0307      \\
  &RMSE       &      -        &    0.0234            &      0.0403       &
     0.0384         &    0.1191          &     0.0487            &
0.1660      \\ \hline
\end{tabular}\\[0.5ex]
\end{center}
\end{table}

Since both the  means and the standard deviations $D(\hat{\bf A},{\bf
A})$ are very small, the estimation for  ${\bf A}$ is accurate.
The coefficients in each CUC models were also estimated accurately.
The errors in estimation decrease as the sample size increases
from 500 to 1000.

Since biases reported in Table~1 are always negative; see also
Figure~1. This indicates that the coefficients in the GARCH(1, 1)
models for CUCs were slightly underestimated. Also note that the
estimation errors decrease when the volatility persistence
(measured by $\alpha_i + \beta_i$) increases; see the upper panel
of Figure~1 for the estimation with the sample size 1000. To make
a comparison, the estimation errors of the GARCH coefficients when
the true ${\bf A}$ is used are plotted in the lower panel.  The
differences are small.

\section{Real data examples}

In this section we illustrate the proposed method
with two real data sets.

The first data set, denoted as SCI, consists of the 2275
daily log returns (in percentages) of S\&P 500 index, stock price
of Cisco System and stock price of Intel Corporation  in  2
January 1991 --- 31 December 1999. This data set has been analyzed
in Tsay~(2001). Figure 2 depicts the time series plots of the
three series.  Descriptive statistics are listed in Table 2.
Obviously, the unconditional distribution of all of these series
exhibit excessive kurtosis; indicating significant departure from
normal distributions.

The Ljung-Box $Q$ statistics suggest some plausible autocorrelation
in these series. But this may be due to the
heteroscedasticity. Hence we compute the $p$-values of these $Q$ tests
based on a bootstrap procedure: for each of the mean-deleted
component return series, we first fit a univariate
GARCH(1,1) model
\begin{equation}\nonumber
Y_t=\sigma_t \epsilon_t, \hspace{1cm} \sigma_t^2=\alpha_0
+ \alpha_1 Y_{t-1}^2 +\beta_1 \sigma_{t-1}^2,
\end{equation}
and denote the estimated
parameters as $\hat{\alpha}_0, \hat{\alpha}_1, \hat{\beta}_1$,
respectively, and the standardized residuals as
$\hat{\epsilon}_t$. Draw $ \epsilon_t^{\ast}$
randomly with replacement from $\{\hat\epsilon_t,\; t=1,\cdots,n\}$ and
draw $Y_t^{\ast}$
from
\begin{equation}\nonumber
Y_t^{\ast}=\sigma_t \epsilon_t^{\ast}, \hspace{1cm}
\sigma_t^2=\hat{\alpha}_0 + \hat{\alpha}_1 Y^{\ast 2}_{t-1}
+\hat{\beta}_1 \sigma_{t-1}^{\ast 2}.
\end{equation}
Let $Q^{\ast}$ be a $Q$-statistic based on $Y_t^{\ast}$.
The $p$-value of $Q$ is now
estimated by the relative frequency of the occurrence of the event
that $Q^{\ast}$ is great than
$Q$ in a repeated bootstrap sampling for 1000 times.
In Table 2,  those $p$-values are listed in parentheses
below the values of the corresponding $Q$ statistics.
Based on those $p$-values, there is no significant evidence
for the existence of autocorrelation in all the three
component series.
Accordingly there is no need to fit a VAR model for the
conditional mean  for this data set.

Let $\bY_t$ be the mean-deleted returns of SCI.
Let $\bSigma =
{\bf P}{\bf \Lambda}{\bf P}^{\tau}$ be the sample covariance
matrix of $\bY_t$, where ${\bf P P}^\tau = \bI_3$ and $\bLambda$
is diagonal. Let ${\bf X}_t={\bf \Lambda}^{-\frac{1}{2}}{\bf
P}^{\tau} {\bf Y}_t$. Then we may regard the (unconditional)
covariance matrix of $\bX_t$ is $\bI_3$.

\begin{table}[htb]
\begin{center}
\caption[Table 2]{Summary Statistics of the Two Real Data Sets }
\begin{tabular}{c | c c c |c c c c c }\hline
               &     S$\&$P 500     &    Cisco          &  Intel            &    HS             &       JN         &   SH              &       ST          &    TW              \\[0.5ex]\hline
N              &       2275         &    2275           &  2275             &   1349            &    1349          &  1349             &      1349         &    1349            \\
Mean           &       0.0656       &   0.2567          &  0.1561           &   -0.0198         &   -0.0477        &  0.0178           &     -0.0081       &  -0.0400           \\
Stdev          &      0.8747        &  2.8540           & 2.4644            &   2.1822          &    1.7382        &  1.5401           &      1.8784       &   1.9863           \\
Min            &      -7.1140       &  -22.1000         & -14.5810          &  -14.7346         &   -9.0145        &  -8.7277          &     -9.1535       &  -9.9360           \\
Max            &       4.9900       &  15.5760          &  12.8500          &   20.2083         &    8.8876        &  8.8491           &     19.5559       &   9.7871           \\
Skewness       &     -0.3600        & -0.3963           &  -0.2353          &    0.6419         &    0.1375        &  0.1861           &      0.9114       &   0.1345           \\
Kurtosis       &      9.0469        &  6.7229           &  5.4701           &   14.3999         &    5.0891        &  8.4310           &     15.2063       &   5.4082           \\ \hline
$Q(10)$        &  22.8322           & 25.3861           & 6.8567            &   32.2251         &    8.8471        & 12.9372           &     28.6943       &  16.9723           \\
               &  \fsize (0.2440)   & \fsize (0.0870)   &  \fsize (0.8180)  & \fsize (0.1760)   &  \fsize (0.7540) &  \fsize (0.7770)  &  \fsize (0.2180)  &   \fsize (0.2540)  \\
$Q(20)$        &  44.2898           & 33.9490           & 30.3427           &   46.1651         &    19.1511       & 26.9255           &     40.7220       &  28.4664           \\
               &  \fsize (0.2300)   & \fsize (0.2500)   & \fsize (0.1170)   & \fsize (0.2810)   &  \fsize (0.7200) & \fsize (0.7310)   &    \fsize (0.2870)&  \fsize (0.3290)   \\ \hline
\end{tabular}\\[0.5ex]
\end{center}
\begin{singlespace}
\emph{Note:} {\sl  $Q(k)$ is referred to the Ljung-Box portmanteau test  statistics.
Figures in parentheses are their corresponding p-values based on 1000 bootstrap
replications.  }
\end{singlespace}
\end{table}

Based on data $\bX_t$, an estimator $\wh \bA $ was obtained with
$\Psi_n(\hat{\bf A})= 0.1732$. Consequently a GARCH(1,1) model was
fitted for each CUC. The estimated coefficients are listed
in Table~3 which shows
that the volatility of the first and third CUCs is highly persistent
as $\hat{\alpha}_1+\hat{\beta}_1=0.9925$ and
$\hat{\alpha}_3+\hat{\beta}_3=0.9611$.
(One may fit the first CUC with an IGARCH model.)
On the other hand, the volatility of the second CUC is less persistent as
$\hat{\alpha}_2+\hat{\beta}_2=0.80$.

We applied the bootstrap procedure (with 500 replications) described
in section~2.4 to test the existence of the CUCs. The $p$-value
is 0.60, indicating that there is no strong evidence against the
hypothesis of the existence of CUCs.
The $(1-\alpha)$  bootstrap
confidence set for the transformation
matrix ${\bf A}$  is
$\{ \bA | D(\hat{\bf A}, {\bf A})\le c_{\alpha}, \;
{\bf A}^{\tau}{\bf A}={\bf I}_3 \}$ with $c_{\alpha}= 0.1718$ for
$\alpha = 0.05$, and 0.1368 for $\alpha = 0.1$.
Since $D(\hat{\bf A}, \bI_3) = 0.2593$,
$\bI_3$ is not contained in the confidence sets. This indicates that the
principal components cannot be taken as the CUCs.
The confidence intervals for the parameters for each CUC-GARCH(1,1)
models are listed in Table~3.
The length of the confidence intervals increase as
the volatility persistent measured by $\wh \alpha_i + \wh \beta_j$
decreases. This is  consistent with the finding from the simulation
study reported  in section 3.

Based on the fitted conditional variances $\wh \sigma_{ti}^2$ for the CUCs,
the conditional variance matrix for the original series $\bY_t$ is equal to
\[
\hat{\bf H}_t={\bf W} \diag\{\wh\sigma^2_{t1}, \wh\sigma^2_{t2},
\wh\sigma^2_{t3}\} {\bf W}^{\tau},
\]
where  ${\bf W}={\bf P}{\bf \Lambda}^{\frac{1}{2}}{\wh\bA}$.
Since the volatility processes of the first and third CUC are highly
 persistent, they can be modelled with Integrated GARCH models. If so,  the volatility processes for original series and their covariance processes are virtually
modelled by  mixtures of IGARCH models and mean-reverting GARCH models,
which is similar to the Component GARCH model used in Ding and
Granger (1996) to capture the long memory properties for a univariate
volatility process.

\begin{table}[htb]
\begin{center}
\caption[Table 4]{ Fitted CUC-GARCH(1,1) model for
SCI }
\begin{tabular}{c | c c c }\hline
               &        Estimate                         &           95\% Confidence Set      &             90\% Confidence Set               \\[0.5ex]\hline
  ${\bf a}_1$  & $(-0.5605,  -0.0018,  -0.8081)^{\tau}$  &                                    &                                                \\
  ${\bf a}_2$  & $(0.5693,   0.7217,   -0.3939)^{\tau}$  &          $c_{0.05}= 0.1718$        &         $c_{0.10}=0.1368$                      \\
  ${\bf a}_3$  & $(0.6015,   -0.6922,  -0.3989)^{\tau}$  &                                    &                                                \\  \hline
  $\gamma_1$   &            0.0074                       &     (0.0042, 0.0592)               &          (0.0048, 0.0449)                      \\
  $\alpha_1$   &            0.0519                       &     (0.0316, 0.0915)               &          (0.0350, 0.0812)                      \\
  $\beta_1$    &            0.9406                       &     (0.8446, 0.9576)               &          (0.8740, 0.9548)                      \\   \hline
  $\gamma_2$   &            0.1997                       &     (0.0460, 0.7138)               &          (0.0673, 0.5705)                      \\
  $\alpha_2$   &            0.0432                       &     (0.0077, 0.1054)               &          (0.0107, 0.0926)                      \\
  $\beta_2$    &            0.7572                       &     (0.2446, 0.9289)               &          (0.3600, 0.9069)                      \\   \hline
  $\gamma_3$   &            0.0389                       &     (0.0200, 0.1042)               &          (0.0239, 0.0870)                      \\
  $\alpha_3$   &            0.0884                       &     (0.0476, 0.1305)               &          (0.0517, 0.1236)                      \\
  $\beta_3$    &            0.8727                       &     (0.7889, 0.9266)               &          (0.8051, 0.9140)                      \\ \hline
\end{tabular}\\[0.5ex]
\end{center}
\end{table}

Figure 3 depicts the fitted volatility processes for each return
series and Figure 4 displays the conditional correlations among
the three components series. Note the volatilities of the S$\&$P
500  index has a much smaller scale than those of the two
individual stocks.
Increasing trends can be observed in all the three correlation processes
over the last three years,  which may be  connected  with the
high volatilities in all the return series over the same period. But on
the other hand, the high volatility of Cisco prices in the middle period did
not lead to a high correlation with the other two. This suggests a unilateral
impact from the market to the single stock.

Figure 5 displays the fitted volatility processes for the three return
series based on the orthogonal GARCH(1,1) model of Alexander (2001) and
Ding and Engle (2001). Note that orthogonal GARCH model effectively
treats the principal
components as conditional uncorrelated variables, which may overlook important
conditional dependence structure in the original data.
Note that the time varying patterns in the three processes in Figure~5 are
similar, which is different from Figure~3 of CUC-GARCH(1,1) fitted.
Especially the orthogonal GARCH fitting artificially inflates
the volatility of S\&P500 index in the middle period; see the original
time plot of the series in Figure~2.
 The inflation is due to treating the conditional correlated principal
components as CUCs. As we stated above, the identity matrix is indeed
not included in the confidence set for~$\bA$.

\askip

Our second data set consists of the daily close returns of five Asian
stock  indices, namely,  Hang Seng index of Hong Kong (HS), Japan Nikkei
225 index (JN),
Shanghai Composite index of China (SC), Straits Time index of Singapore
(ST) and Taiwan Weighted index (TW) in the period of 1 August 1997 ---
30 December 2003.  Adjustments are also made to
account for the differences in the holidays of the five markets.
The five return series are plotted in Figure~6, and the descriptive
statistics are listed in Table~2.
Most of the sample means of these returns are
negative, except the mean of SC. Different from the three series
in SCI, all
five series are right-skewed over this specific
period. The bootstrap $p$-values for the $Q$ statistics are obtained in the
same way as before; indicating no significant
autocorrelation in all the five series.

We fitted a CUC-extended GARCH(1,1) to the mean-deleted return series.
The lagged valued from the other CUCs were selected using BIC together
with a forward searching; see section~2.3.3.
The fitted extended GARCH(1,1) models, based on quasi-MLE with Gaussian
likelihood, for the five
CUCs are reported in Table~4. According to the fitted models, the first
CUC is causal in variance to the fifth CUC, the second CUC is causal in
variance to the first
and the third CUCs, and the fifth CUC is causal in variance to the first CUC.
On the other hand, no additional variables were selected in the models
for the second and fourth CUCs.

Figure~7 displays the fitted volatility processes for the five
original stock returns. As expected, the most volatile waves are
observed at the early of 1998 with the onset of the Asian
financial crisis, which are especially predominant in Hong Kong
and Singapore markets. While the shock is still big, the impact of
the crisis on Japan and Taiwan markets is less drastic.
Furthermore, the effect to Shanghai market is on a much smaller
scale. In Figure 8, we present the fitted conditional correlation
between Hong Kong and the other four markets. Obviously, the most
correlated period is in accord with the epidemic of Asian
financial crisis. After that, the correlations between Hong Kong
and Singapore almost remain at a constant level except two
downslides in the middle of 1999 and 2002, respectively. Likewise,
the correlations between Hong Kong and Taiwan are almost at a
constant level, although a little smaller than
 that with Singapore market.  A upward trend can be seen in the
correlation between Hong Kong and Japan markets in the last few years,
which suggests that these
  two markets  were  becoming more closely integrated.  On the contrary,
the correlations between Hong Kong and Shanghai markets seems to have a
downward to zero trend in
  the last few years. The implications of these observations  to
international diversification deserve a further investigation.

\begin{table}[htb]
\begin{center}
\caption[Table 4]{Extended GARCH(1,1) for CUCs of Asian
Market Data }
\begin{tabular}{c | c | c | l |c  }\hline
     $j$       &   $j_i$     &   $ r$ &        \multicolumn{1}{c|}{                                  $ \sigma_{t,j}  $       }                  &       BIC       \\[0.5ex]\hline
       1       &   5, 2      &    2   &   $\sigma_{t,1}^2=0.0271+0.8609\sigma_{t-1,1}^2+0.0405Z_{t-1,1}^2+0.0637Z_{t-1,5}^2+0.0117Z_{t-1,2}^2 $&    $3622$       \\
       2       &             &    0   &   $\sigma_{t,2}^2=0.0521+0.8004\sigma_{t-1,2}^2+0.1475Z_{t-1,2}^2 $                                    &    $3602$       \\
       3       &    2        &    1   &   $\sigma_{t,3}^2=0.0077+0.9301\sigma_{t-1,3}^2+0.0526Z_{t-1,3}^2+0.0098Z_{t-1,2}^2$                   &    $3731$       \\
       4       &             &    0   &   $\sigma_{t,4}^2=0.0704+0.8539\sigma_{t-1,4}^2+0.0757Z_{t-1,4}^2 $                                    &    $3780$       \\
       5       &    1        &    1   &   $\sigma_{t,5}^2=0.0122+0.8227\sigma_{t-1,5}^2+0.1530Z_{t-1,5}^2+0.0261Z_{t-1,1}^2$                   &    $2534$       \\  \hline
 \end{tabular}\\[0.5ex]
\end{center}
\end{table}

Finally we compared the fitting based on our CUC-based GARCH(1,1) with
the orthogonal GARCH(1,1) models and Engle's dynamic conditional correction
(DCC) model (\ref{a4}) and
(\ref{a5}) in terms of a goodness-of-fit tests based on the Ljung-Box statistic
(Tse and Tsui 1999). Note the DCC-model for each component of $\bY_t$
reduced to the standard univariate GARCH(1,1) fitting.
We define the standardized residual for the
$i$-th series as
$
\hat{u}_{ti}=Y_{ti}/\hat{\sigma}_{t,ii}^{1/2},
$
where $\hat{\sigma}_{t,ii}$ is the $(i,i)$-th element of the fitted
conditional variance of $\bY_{t}$. Define
\[
C_{t,ij}=\left\{ \begin{array} {l r} \hat{u}_{ti}^2-1 & i= j \\
\hat{u}_{ti}\hat{u}_{tj}-\hat{\rho}_{t,ij}  & i\ne j, \end{array} \right.
\]
where
$\hat{\rho}_{t,ij}=\hat{\sigma}_{t,ij}/(\hat{\sigma}_{t,ii}
\hat{\sigma}_{t,jj})^{1/2}$
is the estimated conditional correlation between $Y_{ti}$ and $Y_{tj}$.
If the model
is correctly specified, there is no autocorrelation in $\{ C_{t,ij}, t\ge 1\}$
for any fixed $i, j$.
Put
\[
Q(ij, M)=n \sum_{k=1}^M r_{ij,k}^2,
\]
where $r_{ij,k}$ is the lag $k$ sample autocorrelation of
$C_{t,ij}$. It is intuitively clear that the large values of $Q(ij, M)$ indicate
the lack of fit for the conditional correlation between the $i$-th and
$j$-th components $\bY_t$ for $i\ne j$, and the lack of fit for the
conditional variance of the $i$-th component for $i=j$.  Although the
distribution theory of $Q(ij,M)$ is
unknown, empirical evidence suggests that $\chi^2_M$ provides a
reasonable reference in practice; see Tse and Tsui (1999).

Table~5 lists the values of the $Q$-statistics with $M=10$. The
significant levels were gauged  according to the
$\chi^2_{10}$-distribution. The advantage  of using the CUC-GARCH
model over the Orthogonal GARCH model is obvious as the $Q$-values
for the former tend to be smaller, or significantly smaller, than
those for the latter. Furthermore, all the $Q$ values for the
fitted CUC-GARCH models are insignificant at the level of 10\%,
while the test rejects some Orthogonal GARCH fittings at the
significance level 1\%. For example, the $p$-values for testing
the correlations between S\&P~500 and Cisco stock,  and S\&P~500
and Intel stock is less than 1\%; indicating significant
autocorrelation. This may explain the incomprehensible jumps in
the fitted volatility for S\&P 500 by orthogonal GARCH model in
Figure~5. The same phenomena may also be observed in the fitting
for the second data set. The  orthogonal GARCH model failed to
provide adequate fittings for Hang Seng index (HS), Singapore
Straits Time index (ST) and Taiwan Weighted index (TW), as
indicated by  the large $Q$-values; see Table~5.

Overall the DCC model provide a competitive performance
to the CUC model for the Asian Markets data. This is may due
to a certain degree of homogeneity
among the five Asian market indices.
For SCI consisting of one market index and two stock prices,
the gain of using CUC over DCC is more pronounced. First, the DCC-model
 seems to fail to catch the dynamic correlation
between the returns of the S\&P 500 index and the Cisco stock price.
Furthermore, although $Q$-value for the CUC-model for S\&P 500 is
marginally larger than that of the DCC model,
the $Q$-values for the CUC-models for both Intel and Cisco prices
are substantially smaller than those for the DCC models; suggesting
an improvement for the modelling  volatility dynamics for the Intel or
the Cisco price by incorporating the information from other series.

The $Q$-tests with different values of $M$ lead to similar pattern
as Table~5, which, therefore, are omitted to save the space.

\begin{landscape}
\begin{table}[tabh]
\begin{center}
\caption[Table 5]{Specification test ---  $Q(10)$ for cross products of standardized residuals }
\begin{tabular}{@{\hspace{0.6cm}}c @{\hspace{0.6cm}} | @{\hspace{0.5cm}} c @{\hspace{0.8cm}} c @{\hspace{0.8cm}} c @{\hspace{0.5cm}} |  c @{\hspace{1cm}} c @{\hspace{1cm}}c @{\hspace{1cm}}c  @{\hspace{0.6cm}} }\hline
               &             \multicolumn{3}{c|} { SCI   data}      &                     \multicolumn{4}{c} {Asian Market Data}                   \\ \hline
   $i,j$       &    O-GARCH       &   DCC          &   CUC-GARCH    &   \hspace{0.5cm} O-GARCH      &    DCC         &    CUC-GARCH    &   CUC-Ex GARCH      \\[0.5ex]\hline
     1         &  $59.9140^{***}$ &   5.9498       &  6.2050        & \hspace{0.5cm}$56.7580^{***}$ &    6.0285      &    11.4480      &   8.6961                  \\
     2         &   10.5100        &   9.0587       &  8.0542        & \hspace{0.5cm}12.3540         &    7.8517      &    8.6713       &   8.7751                  \\
     3         &   2.6192         &   6.4293       &  2.2397        & \hspace{0.5cm} 8.5368         &    9.2749      &    8.5301       &   8.5265                  \\
     4         &                  &                &                & \hspace{0.5cm}$18.6100^{**}$  &    2.6610      &    4.0512       &  3.7954                   \\
     5         &                  &                &                & \hspace{0.5cm}$18.0610^{*}$   &    7.4710      &    11.7960      &  13.5150                  \\
     1,2       &  $51.8060^{***}$ &  10.4887       &  10.9090       & \hspace{0.5cm} 7.1025         &    7.0622      &    4.6433       &   4.2671                  \\
     1,3       &  $77.5140^{***}$ & $20.6745^{**}$ &  10.5170       & \hspace{0.5cm} 3.8940         &    4.6465      &    3.4987       &   3.5564                  \\
     1,4       &                  &                &                & \hspace{0.5cm}$17.2180^{*}$   &    4.7943      &    6.2915       &   5.8084                  \\
     1,5       &                  &                &                & \hspace{0.5cm} 9.2396         &    6.1648      &    5.6669       &   6.3143                  \\
     2,3       &   5.9453         &   7.0617       &  9.6275        & \hspace{0.5cm} 9.6031         &    10.1762     &    9.6444       &   9.5912                  \\
     2,4       &                  &                &                & \hspace{0.5cm} 6.3708         &    7.7241      &    3.4542       &   3.2648                  \\
     2,5       &                  &                &                & \hspace{0.5cm} 6.8629         &    5.8438      &    6.1856       &   6.9089                  \\
     3,4       &                  &                &                & \hspace{0.5cm} 11.9120        &    8.0303      &    7.3119       &   5.8486                  \\
     3,5       &                  &                &                & \hspace{0.5cm} 2.2256         &    2.1565      &    1.5721       &   1.6857                  \\
     4,5       &                  &                &                & \hspace{0.5cm} 5.4389         &    4.7838      &    3.0312       &   3.1083                  \\  \hline
 \end{tabular}\\[0.5ex]
\end{center}
\begin{singlespace}
\emph{Note:} {\sl 1)  ***, **, * indicate that the corresponding
test is significant at the level 0.01, 0.05 and 0.1, respectively.
2)  $i, j$ in the left column corresponds to to the orders of component
series in each data sets. For example, ``1,2'' stands for the cross
product of the standardized residuals of S\&P 500 and Cisco for the SCI,
and for HS and JN for the Asian market data set.}
\end{singlespace}
\end{table}
\end{landscape}

 \setcounter{equation}{0}
 \renewcommand{\theequation}{5.\arabic{equation}}

\section*{Appendix A --- Proof of Theorem 1}

We introduce some notation first.
Let $$\bC_{n, k}(B) = (n-k)^{-1} \sum_{t=k+1}^n \bX_t \bX_t^\tau I
(\bX_{t-k} \in B), \quad \quad \bC_k (B) = E\{ \bX_t \bX_t^\tau I (\bX_{t-k} \in
B)\}.$$
The lemma below shows that both $\Psi(\cdot)$ and $\Psi_n(\cdot)$ are
Lipschitz continuous on $\calH_D$ with $D$-distance, where $\calH_D$ is
the quotient space; see Remark 2.

\askip

\noindent
{\bf Lemma 1}. For any $\bU, \bV \in \calH_D$, it holds that
 $$|\Psi(\bU) - \Psi(\bV) | \le c \; \tr E (\bX_t \bX_t^T) \, \{ D(\bU, \bV) \}^{1/2},$$ and
$$
|\Psi_n(\bU) - \Psi_n(\bV) | \le c \; \tr (n^{-1} \sum_{i=1}^n
\bX_t \bX_t^T ) \, \{ D(\bU, \bV) \}^{1/2}$$ almost surely, where
$c>0$ is a constant and $\tr(\bA)$ is the trace of a matrix $\bA$.

\askip

\noindent {\bf Proof}. We only prove the lemma for $\Psi(\cdot)$.
The result for $\Psi_n(\cdot)$ may be shown in the same manner.
Let $\bU=(\bu_1, \cdots, \bu_d)^\tau$, $\bV=(\bv_1, \cdots,
\bv_d)^\tau$, $u_{ijk}(B) = E\{ \bu_i^\tau \bC_k(B) \bu_j\}$ and
$v_{ijk}(B) = E\{ \bv_i^\tau \bC_k(B) \bv_j\}$. We assume that the
orders and the directions of $\bu_i$ and $\bv_j$ are arranged such
that $\bu_i ^\tau \bv_i\in [0,1]$ for all $i$, and
\begin{equation} \label{p1}
D(\bU, \bV) = 1 - {1\over d} \sum_{i=1}^d \bu_i ^\tau \bv_i
= {1\over d} \sum_{i=1}^d(1 - \bu_i ^\tau \bv_i).
\end{equation}
See (\ref{b7}).
Put the spectral decomposition for $\bC_k(B)$ as
$$\bC_k(B) = \sum_{\ell=1}^d \mu_{\ell}(B,k) \bgamma_\ell
\bgamma_\ell^\tau,$$ where $\mu_1(B,k) \ge \cdots \ge
\mu_d(B,k)\ge 0$ are the eigenvalues of $\bC_k(B)$, and
$\bgamma_1, \cdots, \bgamma_d$ are their corresponding (orthonormal)
eigenvectors. It is easy to see that $\mu_\ell(B,k)\le \mu_\ell$
for all $k$ and $B$, where $\mu_1 \geq \cdots \geq \mu_d$ are the
eigenvalues of the matrix $E \{ \bX_t \bX^\tau\}$.
Consequently, by noticing that
$|\bgamma_\ell^\tau \bu_j | \leq 1$ and $|\bv_i^\tau \bgamma_\ell |
\leq 1$, we have
\begin{eqnarray} \nonumber
&& | u_{ijk}(B) - v_{ijk}(B) | \; \le \; \sum_{\ell=1}^d
\mu_\ell | \bu_i^\tau \bgamma_\ell \bgamma^\tau_\ell \bu_j -
\bv_i^\tau \bgamma_\ell \bgamma^\tau_\ell \bv_j|\\ \nonumber
&\le & \sum_{\ell=1}^d
\mu_\ell \{ | \bu_i^\tau \bgamma_\ell \bgamma^\tau_\ell \bu_j -
 \bv_i^\tau \bgamma_\ell \bgamma^\tau_\ell \bu_j|
+ | \bv_i^\tau \bgamma_\ell \bgamma^\tau_\ell \bu_j -\bv_i^\tau
   \bgamma_\ell \bgamma^\tau_\ell \bv_j|\}\\ \nonumber
&\le &
\sum_{\ell=1}^d \mu_\ell \{ | ( \bu_i-\bv_i)^\tau \bgamma_\ell| +
|\bgamma^\tau_\ell (\bu_j-\bv_j)|\}
\end{eqnarray}
By using the Cauchy-Schwartz's inequality, the above inequality is
furthered bounded by
\begin{eqnarray}
& & \sum_{\ell=1}^d \mu_\ell  \{ ||\bu_i-\bv_i|| +
||\bu_j-\bv_j||\} \nonumber \\ \label{p2} &=& \sqrt{2} \{( 1 -
\bu_i^\tau \bv_i)^{1/2} + (1 - \bu_j^\tau \bv_j)^{1/2}\}
\sum_{\ell=1}^d \mu_\ell.
\end{eqnarray}

Note that for $x\ne 0$, it holds that
\begin{equation} \label{p3}
|x+y| - |x| = y\, \sgn(x) + 2(x+y)\{ I(-y<x<0) - I(0<x<-y) \}.
\end{equation}
Hence,
\begin{eqnarray}
&&
\Psi(\bU) \; = \;
\sum_{1\le i < j \le d}
\sup_{1\le k \le k_0, \, B\in \calB} \big[
|v_{ijk}(B)|+|v_{ijk}(B) + \{u_{ijk}(B)-v_{ijk}(B)\}|
- | v_{ijk}(B)| \big] \nonumber \\
&=& \sum_{1\le i < j \le d} \sup_{1\le k \le k_0, \, B\in \calB}
\big[ | v_{ijk}(B)| + \{ u_{ijk}(B)-v_{ijk}(B)\}\sgn\{
v_{ijk}(B)\}  \nonumber \\
& & + \; 2 u_{ijk}(B) \{I (B_1) - I( B_2) \}\big], \label{p4}
\end{eqnarray}
where
\[
B_{1} = \{ v_{ijk}(B)-u_{ijk}(B) < v_{ijk}(B) <0\}, \quad B_{2} =
\{ 0< v_{ijk}(B)<v_{ijk}(B)-u_{ijk}(B)\} .
\]
On the set $B_1 \cup B_2$,
\[
|u_{ijk}(B)| \le |u_{ijk}(B)-v_{ijk}(B)|+|v_{ijk}(B)|\le 2
|u_{ijk}(B)-v_{ijk}(B)|.
\]
This, combining with (\ref{p2}) and (\ref{p4}), implies that
\begin{eqnarray}\nonumber
& & |\Psi(\bU) - \Psi(\bV) | \\
&  \le & \sum_{1\le i < j \le d} \sup_{1\le k \le k_0, \, B\in
\calB} \big[\sqrt{2} \{ ( 1 - \bu_i^\tau \bv_i)^{1/2}  + (1 -
\bu_j^\tau \bv_j)^{1/2} \} \sum_{\ell=1}^d \mu_\ell
     +  2|u_{ijk}(B)| I_1 (B_1) \big] \nonumber \\
& \le & \; 5 \sqrt{2} \sum_{1\le i < j \le d} \{( 1 - \bu_i^\tau
\bv_i)^{1/2} + (1 - \bu_j^\tau \bv_j)^{1/2}\} \sum_{\ell=1}^d
\mu_\ell \nonumber \\
& \le &  10 \sqrt{2} d \sum_{\ell=1}^d \mu_\ell \sum_{i=1}^d (1-
\bu_i^\tau \bv_i)^{1/2}. \label{p5}
\end{eqnarray}

Now the lemma follows from (\ref{p5}) and the inequality
\[
\sum_{i=1}^d (1- \bu_i^\tau \bv_i) ^{1/2} \le d^{1/2} \Big\{
\sum_{i=1}^d (1- \bu_i^\tau \bv_i) \Big\}^{1/2},
\]
see also (\ref{p1}). This completes the proof.

\askip

\noindent
{\bf Proof of Theorem 1}.
Since $\bC_{n, k} (B) - \bC_k(B)$ is a real symmetric matrix, it holds
for any unit vectors $\ba$ and $\bb$ that
\[
|\ba^\tau \{ \bC_{n, k} (B) - \bC_k(B)\} \bb| \le || \bC_{n, k} (B) - \bC_k(B)||,
\]
where $|| \bC_{n, k} (B) - \bC_k(B)||$ denotes the sum of the absolute values of
the eigenvalues of $\bC_{n, k} (B) - \bC_k(B)$. This may be obtained by
using the spectral decomposition of $\bC_{n, k} (B) - \bC_k(B)$.
Consequently it holds uniformly for any orthogonal matrix $\bA$ that
\begin{eqnarray} \nonumber
|\Psi_n(\bA) - \Psi(\bA)| & \leq & \sum_{1 \leq i < j \leq d}
     \sup_{1 \leq k \leq k_0, B \in \calB} \left | \ba_i ^\tau \{ \bC_{n,
     k} (B) - \bC_k(B) \} \ba_j  \right | \\
     & \leq & {d(d-1)\over 2}
     \sup_{1 \leq k \leq k_0, B \in \calB}  \| \bC_{n, k}
     (B) - \bC_k(B) \| .
\label{p6}
\end{eqnarray}
Note the $(i,j)$-th element of $\bC_{n, k} (B) - \bC_k(B)$
is $${1 \over n-k}  \sum_{t=k+1}^n X_{ti}X_{tj}
     I(\bX_{t-k} \in B) - E\{ X_{ti}X_{tj} I(\bX_{t-k} \in B)\},$$
where $X_{ti}$ denotes the $i$-th element of $\bX_t$.
Since $E | X_{ti}X_{tj}| < \infty$ and $\calB$ is a VC-class, the
covering number for the set of functions $\{X_{ti}X_{tj}
I(\bX_{t-k} \in B), B \in \calB\}$ has a polynomial rate of growth
for any underlying probability measure (Theorem 2.6.4, van der
Vaart and Wellner 1996).   Hence, it is a Glivenko-Cantelli class.
It follows now from Theorem 3.4 of Yu (1994) that
\[
    \sup_{B \in \calB} \Big | {1 \over n-k} \sum_{t=k+1}^n
    X_{ti}X_{tj} I(\bX_{t-k} \in B) - E\{ X_{ti}X_{tj} I(\bX_{t-k} \in B)\} \Big |
\scon 0,
\]
Consequently, $$\sup_{B \in \calB} |\lambda_{\max}(B, k)| \scon 0,
\quad
\quad
\sup_{B \in \calB}| \lambda_{\min}(B,k)| \scon 0,$$
where $\lambda_{\max}(B, k)$ and $ \lambda_{\min}(B,k)$ denote, respectively,
the maximum and the minimum eigenvalues of $\bC_{n, k} (B) - \bC_k(B)$.
Thus
$$
    \sup_{B \in \calB}  \| \bC_{n, k} (B) - \bC_k(B) \| \scon 0,
$$
for $k=1, \cdots, k_0$. Now it follows from (\ref{p6}) that
$$
  \sup_{\bA \in \calH_D} | \Psi_n(\bA) - \Psi(\bA) | \scon 0.
$$
Combining this with Lemma~1 above and
the continuity of the argmax mapping (Theorem 3.2.2 and Corollary
3.2.3, van der Vaart and Wellner, 1996),  it
holds that $D(\hat{\bA}, \bA_0)
\scon 0$.  This completes the proof of the first part of Theorem 1.

\askip

Under the additional condition $ E | X_{ti}X_{tj}|^{2p} < \infty$
and the mixing condition given in Condition (A4),  Theorem 1 of
Arcones and Yu (1994) implies that the set of functions
$\{X_{ti}X_{tj} I(\bX_{t-k} \in B), B \in \calB\}$ is a Donsker
class, and hence the process $ \{\bDelta_{n, k}(B), B \in \calB
\}$ indexed by $B \in \calB$  converges weakly to a Gaussian
process, where $\bDelta_{n, k} (B) = \sqrt{n} \{ \bC_{n,k}(B) -
\bC_k(B) \}$. It follows from (\ref{p3}) that
\begin{eqnarray} \nonumber
\Psi_n(\bA) &=& \sum_{1 \leq i < j \leq d} \sup_{B \in \calB, 1 \leq
k \leq k_0 } \big [ |\ba_i^T \bC_k(B) \ba_j| + n^{-1/2}
\sgn\{\ba_i^\tau \bC_k(B) \ba_j\} \ba_i ^\tau \bDelta_{n, k}(B) \ba_j \\
& & + \ba_i^\tau \bC_{n,k}(B) \ba_j \{ I(B_3) - I(B_4)\}  \big ]  \nonumber \\
& = &  \Psi (\bA) + O_P(n^{-1/2}), \label{p7}
\end{eqnarray}
where
\[
B_3 = \{n^{-1/2} \ba_i ^\tau \bDelta_{n, k}(B) \ba_j  < \ba_i^\tau
\bC_k(B) \ba_j<0\}, \quad
B_4=  \{0< \ba_i^\tau \bC_k(B) \ba_j < n^{-1/2}
\ba_i ^\tau \bDelta_{n, k}(B) \ba_j\}.
\]
The last equality in (\ref{p7}) follows from the fact that on
$B_3 \cup B_4$,
\[
|\ba_i^\tau \bC_{n,k}(B) \ba_j| \le |\ba_i^\tau \bC_k(B) \ba_j|
+ n^{-1/2}|\ba_i ^\tau \bDelta_{n, k}(B) \ba_j|
\le 2 n^{-1/2}|\ba_i ^\tau \bDelta_{n, k}(B) \ba_j|.
\]
It follows from (\ref{p7}) and condition (A5) that
\begin{eqnarray}
\Psi_n (\bA_0) - \Psi_n(\bA)   =  \Psi(\bA_0) -
    \Psi(\bA) + O_P(n^{-1/2})
 \leq  -a D(\bA_0, \bA) + O_P(n^{-1/2}) \label{p8}.
\end{eqnarray}
Now by substituting $\bA$ by $\hat{\bA}$, the left hand side of
(\ref{p8}) must be non-negative by the definition of~$\hat{\bA}$.
The right hand side of (\ref{p8}) would be negative unless
$$
D(\bA_0, \hat{\bA}) = O_P(n^{-1/2}).
$$
This completes the proof.

\section*{Appendix B --- Proof of Theorem 2}

From the proof of Theorem 1, we have
\begin{equation}
      \sup_{\bA \in \calH} | \Psi_n (\bA) - \Psi(\bA) | \scon 0.
      \label{p9}
\end{equation}
Since $\Psi(\bA)$ is continuous on the compact quotient space
$\calH$, there exists a minimizer $\bA_0$.  It follows that
\begin{eqnarray*}
 \Psi(\hat{\bA}) - \Psi(\hat{\bB}) & = & \Psi(\bA_0) -
 \Psi(\hat{\bB})  + \Psi(\hat{\bA})- \Psi(\bA_0) \\
 & \leq & \Psi(\hat{\bA})- \Psi(\bA_0) \\
 & = & \{\Psi(\hat{\bA}) - \Psi_n( \hat{\bA})\}
        + \{\Psi_n( \hat{\bA}) - \Psi_n( \bA_0 )\}
        + \{ \Psi_n( \bA_0 ) - \Psi( \bA_0)\}.
\end{eqnarray*}
Using the fact $\Psi_n(\hat{A}) - \Psi_n (\bA_0) \leq 0$, we
conclude from (\ref{p9}) that
$$
\liminf \{ \Psi(\hat{\bA}) - \Psi(\hat{\bB})\} \leq 0.
$$
This completes the proof of Theorem 2.

\section*{Appendix C --- Proof of Theorem 3}

For each $j$, there are at most $r$ non-zero $\alpha_{jk}$. Since
$\beta_j < 1$, it holds that
\[
\sigma_{tj}^2 = {\ga_j \over 1- \beta_j} +  \sum_{i=1}^d
\alpha_{ji}\sum_{k=1}^\infty \beta_j^{k-1} Z_{t-k,i}^2.
\]
Now Theorem~2 follows from Lemma~2 below immediately by letting
$Y_{tj} = X_{tj}^2$ and $\rho_{tj} = \sigma_{tj}^2$. Note that
Lemma~2 may be proved in the similar manner to the proof of
Theorem~1 of Giraitis at al~(2000); see also section~2.7.1 of Fan
and Yao~(2003).

\bigskip

\noindent {\bf Lemma 2}. Consider a $d$-dimensional ARCH($\infty$)
process $\bY_t = (Y_{t1}, \cdots, Y_{td})^\tau$ defined by
\[
Y_{tj} = \rho_{tj} \zeta_{tj}, \quad \quad \rho_{tj} = c_j +
\sum_{i=1}^d \sum_{k=1}^\infty b_{jik} Y_{t-k, i}
\]
for $j=1, \cdots, d$, where $\{ \zeta_{tj} \}$ is a sequence of
non-negative i.i.d. random variables with $E(\zeta_{tj}) =1$,
$Y_{tj}\ge 0$, $c_j, b_{jik} \ge 0$. Furthermore, for each $j$,
$b_{jik} \ne 0$ for at most $r (\ge 0)$ values of $k$. Then the
above model admits a unique strictly stationary solution $\{ \bY_t
\}$ with the finite mean
\[
E(\bY_t) = (\bI_d - \bB)^{-1} (c_1, \cdots, c_d)^\tau
\]
under the condition $
  \max_{ 1\le j,\, i \le d}  b_{ji\,\cdot} < 1/r,
$ where $b_{ji\,\cdot} = \sum_{k\ge 1} b_{jik}$, and $\bB$ is a
$d\times d$ matrix with $b_{ji\,\cdot}$ as its $(j,i)$-th element.

\section*{References}
\begin{description}
\begin{singlespace}
\item Alexander, C. (2001). Orthogonal GARCH. In {\sl Mastering
Risk}. Financial Times-Prentice Hall: London; {\bf 2}, 21-38.

\item Arcones, M.A. and Yu, B. (1994).  Central limit theorems for
       empirical processes and U-processes of stationary mixing sequences.
        {\em Jour. Theor. Probab.}, {bf 7}, 47--71.

\item Back, A. and Weigend, A.S. (1997). A first application on
independent component analysis to extracting structure from stock
returns. {\sl International Journal of Neural Systems}, {\bf 8},473-484.
\item
Bauwens, L., Laurent, S. and Rombouts, J.V.K. (2003). Multivariate GARCH models:
 a survey. {\sl A preprint}.
\item
Bollerslev, T. (1990).  Modelling the coherence in short-run nominal exchange rates: a multivariate generalized ARCH model. {\sl Review of Economics and Statistics},
 {\bf 72}, 498-505.
\item
Bollerslev, T.R., Engle, R. and Wooldridge, J. (1998). A capital asset pricing
model with time varying covariances. {\sl Journal of Political Economy}, {\bf 96},
116-131.
\item Chen, M. and An, H. (1998). A note on the stationarity and the existence of moments
      of the GARCH models. {\sl Statistica Sinica}, {\bf 8}, 505-510.
\item Chow, Y.S. and Teicher, H. (1997). {\sl Probability  Theory} (3rd
edition). Springer, New York.
\item Ding, Z. and Granger, C.W.J. (1996). Modeling volatility persistence of speculative returns:
      A new approach. {\sl Journal of Econometrics}, {\bf 73}, 185-215.
\item Ding, Z. and Engle, R. (2001). Large scale conditional covariance matrix modeling, estimation and testing.
      {\sl  Working Paper}, {\bf FIN-01-029}, NYU Stern School of Business.
\item Engle, R. (2002).Dynamic conditional correlation -- a simple class of multivariate
      GARCH models. {\sl Journal of Business and Economic Statistics}, {\bf 20}, 339-350.
\item
Engle, R.F., Ito, T. and Lin, W.-L. (1990). Meteor shoers or
                heat waves? heteroskedastic intra-daily volatility in
                the foreign exchange market.
                {\sl Econometrica}, {\bf 58}, 525-542.
\item Engle, R.F. and Kroner, K.F. (1995). Multivariate simultaneous generalised ARCH.
      {\sl Econometric Theory }, {\bf 11}, 122-150.
\item Engle, R.F., Ng, V.K. and Rothschild, M. (1990). Asset pricing with a factor ARCH covariance structure:
      Empirical estimates for treasury bills. {\sl Journal of
Econometrics}, {\bf 45}, 213-238.
\item Engle, R.F. and Sheppard, K. (2001). Theoretical and empirical properties
of dynamic conditional correlation multivariate GARCH. {\sl A preprint}.
\item Fan, J. and Yao, Q. (2003). {\sl Nonlinear Time Series: Nonparametric and Parametric Methods}.
      Springer, New York.
\item Giraitis, L., Kokoszka, P., and Leipus, R. (2000). Stationary ARCH models: Dependence structure and
      central limit theorem. {\sl Econometric Theory}, {\bf 16}, 3--22.
\item Hall, P. and Yao, Q. (2003). Inference for ARCH and GARCH models. {\sl Econometrica},
      {\bf 71}, 285-317.
\item Harvey, A., Ruiz, E. and Shephard, N. (1994). Multivariate stochastic
variance models. {\sl The Review of Economic Studies}, {\bf 61}, 247-264.
\item Hyv\"arinen, A., Karhunen, J. and Oja, E. (2001). {\sl Independent Component
      Analysis}. Wiley, New York.
\item Jerez, M., Casals, J. and Sotoca, S. (2001). The likelihood of multivariate GARCH models is ill-conditioned. {\sl A preprint}.
\item Kiviluoto, K. and Oja, E. (1998). Independent component analysis for parallel financial time series.
      In {\sl Proc. Int. Conf. on Neural Information Processing (ICONIP'98)}, vol.2, pp.895-989, Tokyo.
\item M${\breve {\rm a}}$l${\breve {\rm a}}$roiu, S., Kiviluoto, K. and
Oja, E. (2000). Time series prediction with independent component analysis. {\sl A
       preprint}.
\item Mikosch, T. and Straumann, D. (2004). Stable limits of martingale transforms with application to the
      estimation of GARCH parameters. {\sl A preprint}.
\item Peng, L. and Yao, Q. (2003). Least absolute deviations estimation for ARCH
      and GARCH models. {\sl Biometrika}, {\bf 90}, 967-975.
\item Penzer, J., Wang, M. and Yao, Q. (2004). Approximating volatilities by asymmetric power GARCH functions. {\sl A preprint}.
\item Tsay, R. (2001). {\sl Analysis of Financial Time Series}. Wiley, New York.
\item Tse, Y. K. and Tsui, A.K.C. (1999). A note on diagnosing multivariate conditional heteroscedasticity models.
     {\sl Journal of Time Series Analysis}, {\bf 20}, 679-691.
\item Vilenkin, N. (1968). Special functions and the theory of group representation, translations of
      mathematical monographs. {\sl American Math. Soc.}, Providence, Rhode Island, 22.

\item van der Vaart, A.W. and  Wellner, J.A. (1996).  Weak Convergence and
         Empirical Processes.  Springer, New York.

\item van der Weide, R. (2002). GO-GARCH: a multivariate generalized orthogonal GARCH model. {\sl Journal of
      Applied Econometrics}, {\bf 17}, 549-564.

\item Wang, M. and Yao, Q. (2005). Modelling multivariate volatilities: an ad hoc
approach. To appear in ``{\sl Contemporary Multivariate Analysis and
Experimental Designs}'' J. Fan,  G. Li \& R. Li (edit.) World Scientific,
Singapore.

\item Yu, B. (1994).  Rates of convergence for empirical processes
      of stationary mixing sequences.   {\sl Ann. Statist.}, {\bf 22},
     94-116.

\end{singlespace}
\end{description}

\begin{figure}[hb]
\centerline{\psfig{figure=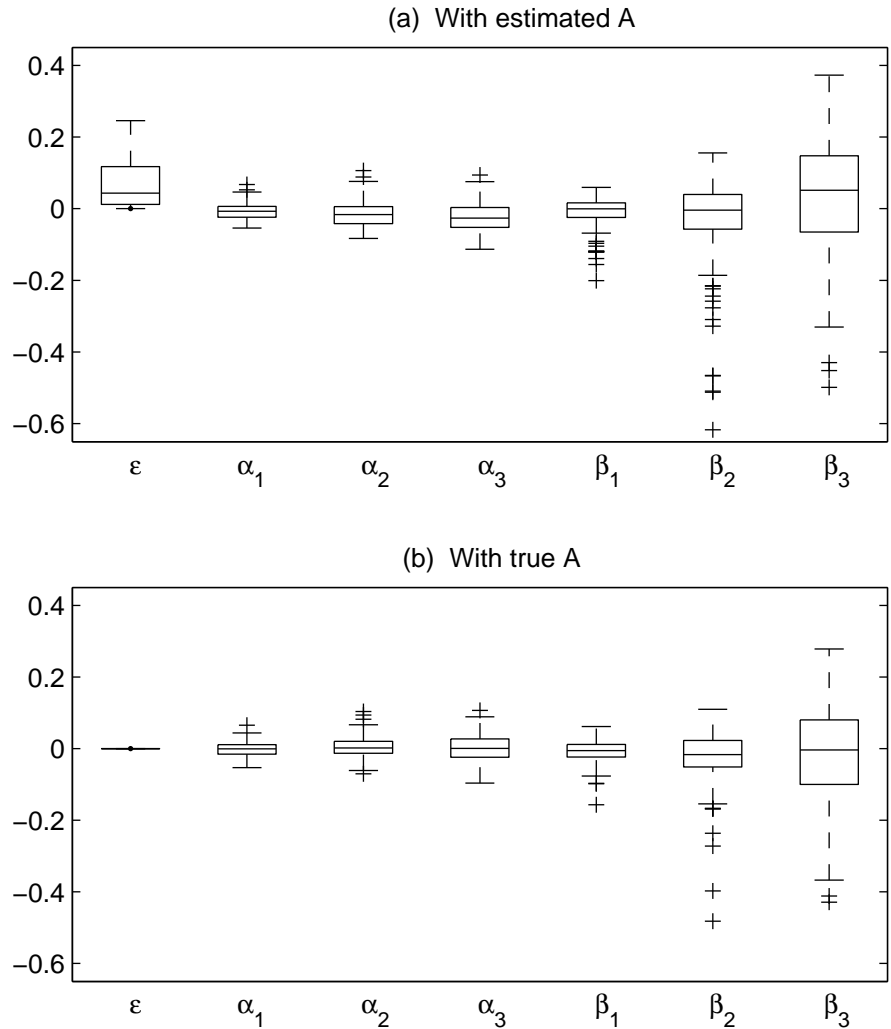}}
\begin{singlespace}

\caption[Fig 1] {\sl Boxplots of the errors in estimation for
CUC-GARCH(1,1) model (\ref{ex1}) with
$\bA =\wh \bA$ estimated (upper panel) and the true $\bA$ (lower panel).
 The sample size is $n=1000$.}
\end{singlespace}
\end{figure}

\newpage

\begin{figure}
\centerline{\psfig{figure=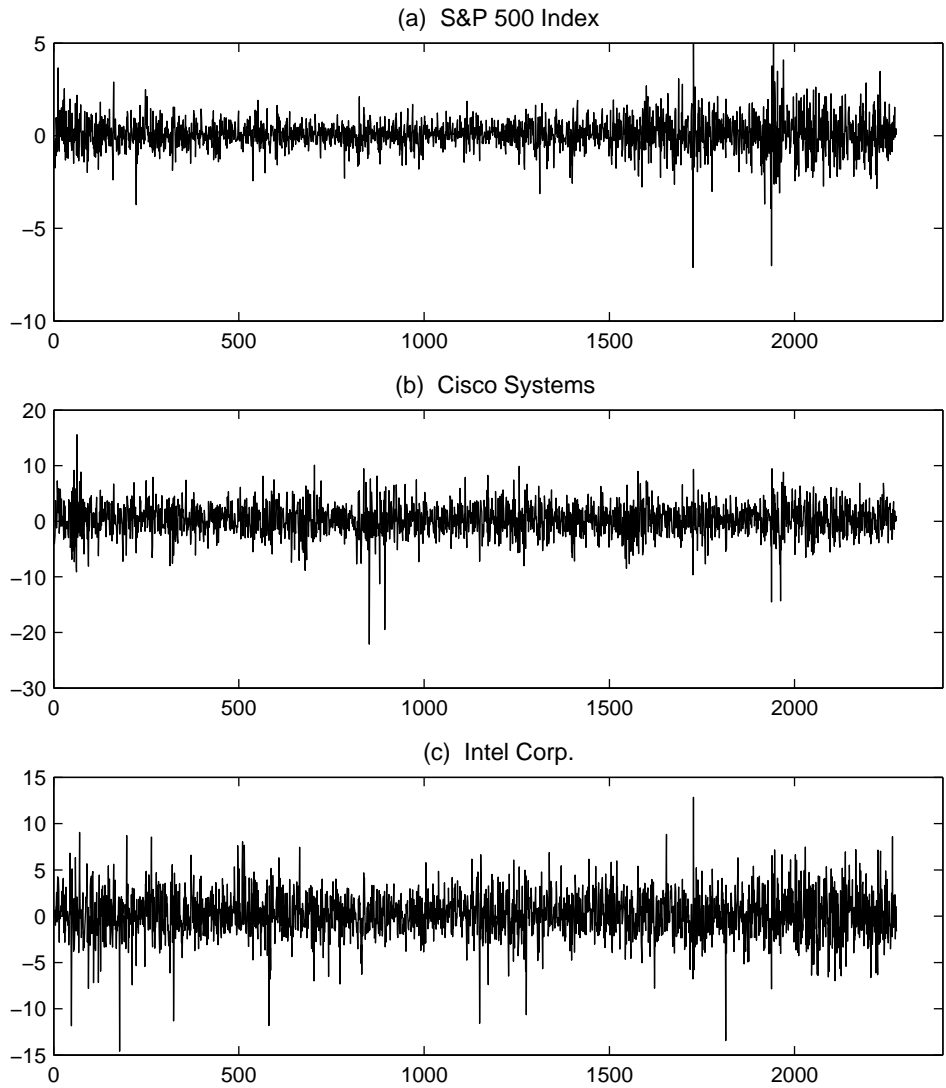}}
\begin{singlespace}

\caption[Fig 2] {\sl Plots of daily log return of
(a)  $S\&P$ 500 index, (b) Cisco Systems stock  and (c) Intel
Corporation stock.  Time span is from January 2, 1991 to December 31, 1999 with 2275
observations.}

\end{singlespace}
\end{figure}

\newpage

\begin{figure}
\centerline{\psfig{figure=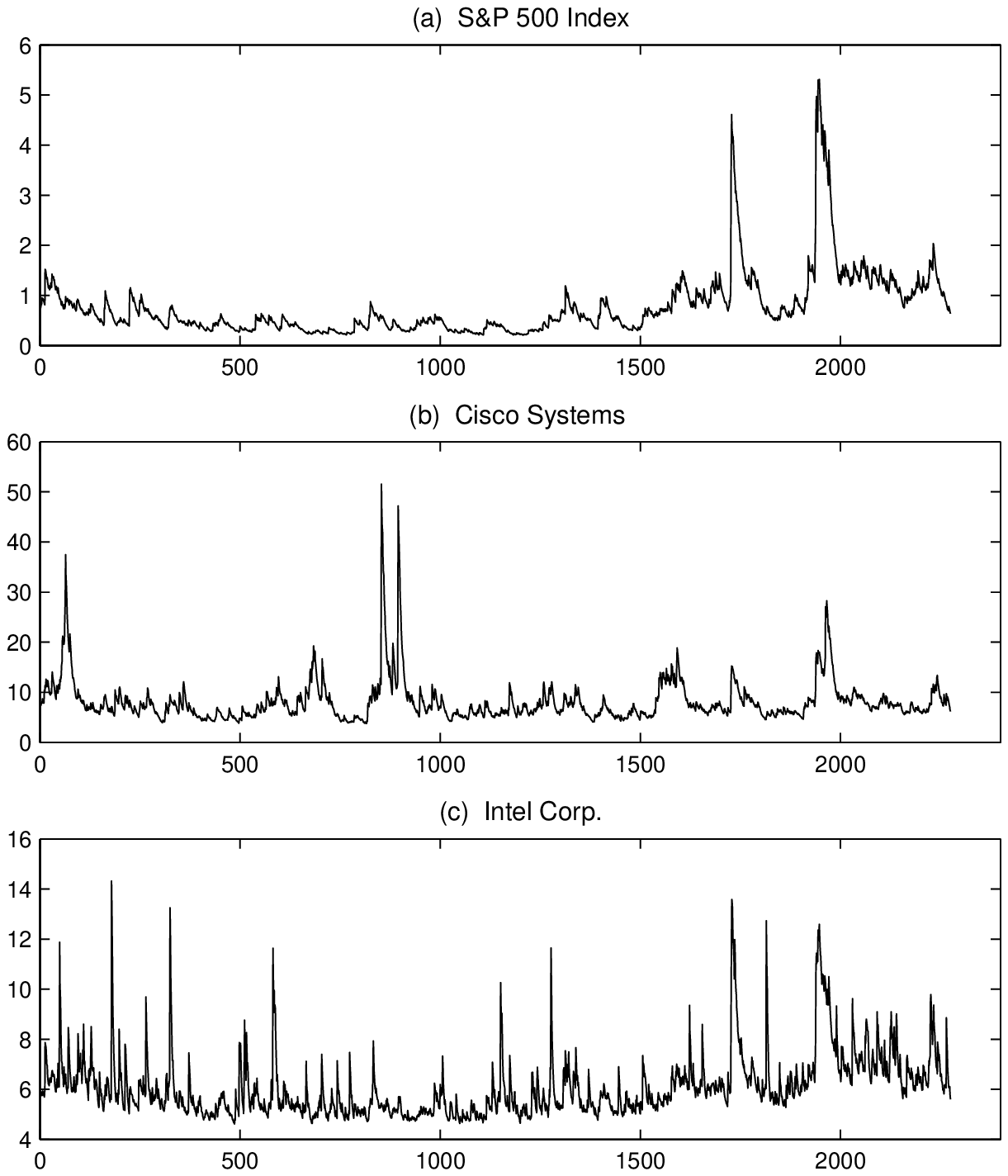}}
\begin{singlespace}

\caption[Fig 3] {\sl Fitted volatility processes based on  CUC-GARCH(1,1) model for daily log returns of
(a)  $S\&P$ 500 index, (b) Cisco Systems stock  and (c) Intel
Corporation stock. }

\end{singlespace}
\end{figure}

\newpage

\begin{figure}
\centerline{\psfig{figure=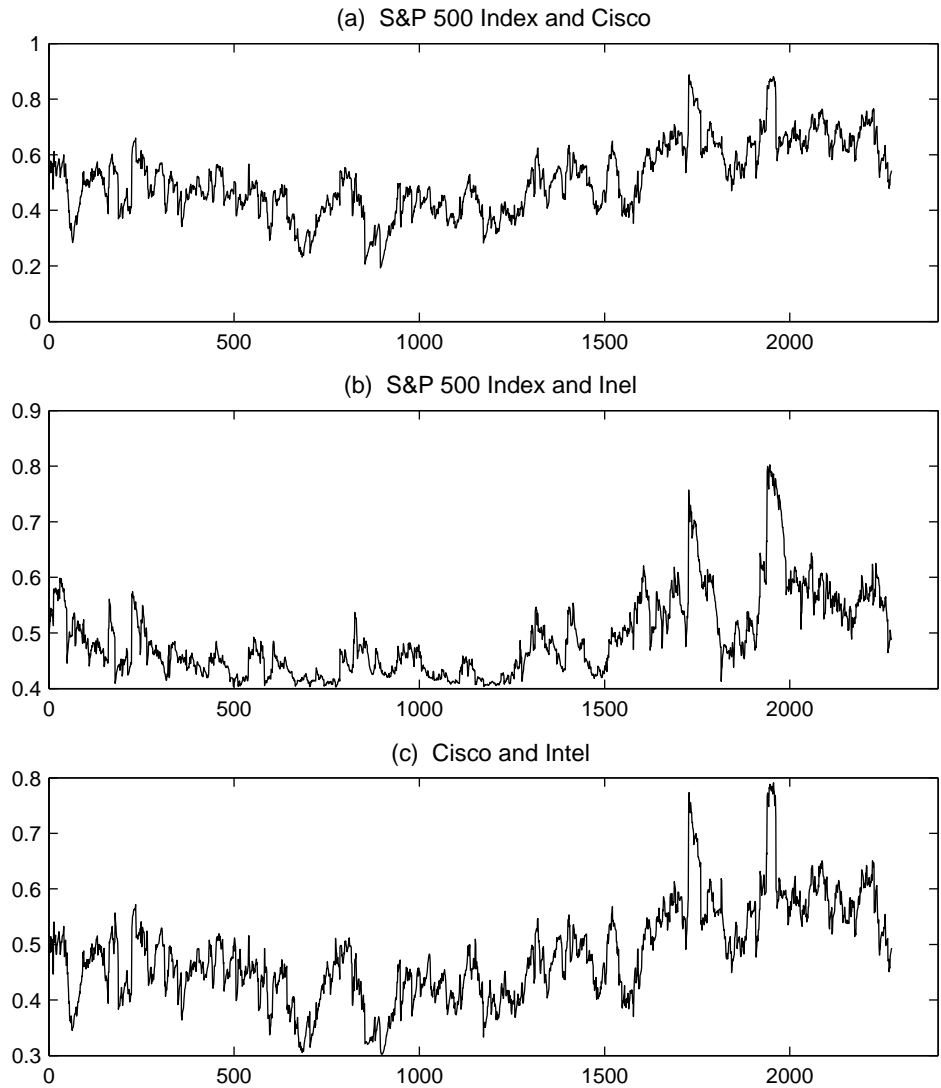}}
\begin{singlespace}

\caption[Fig 4] {\sl Fitted conditional correlations based on  CUC-GARCH(1,1) model
for  daily log returns between (a) $S\&P$ 500 index and Cisco Systems
stock, (b) $S\&P$ 500 index and Intel Corporation stock,  and (c)
Cisco Systems stock and Intel Corporation stock.  }
\end{singlespace}
\end{figure}

\newpage

\begin{figure}
\centerline{\psfig{figure=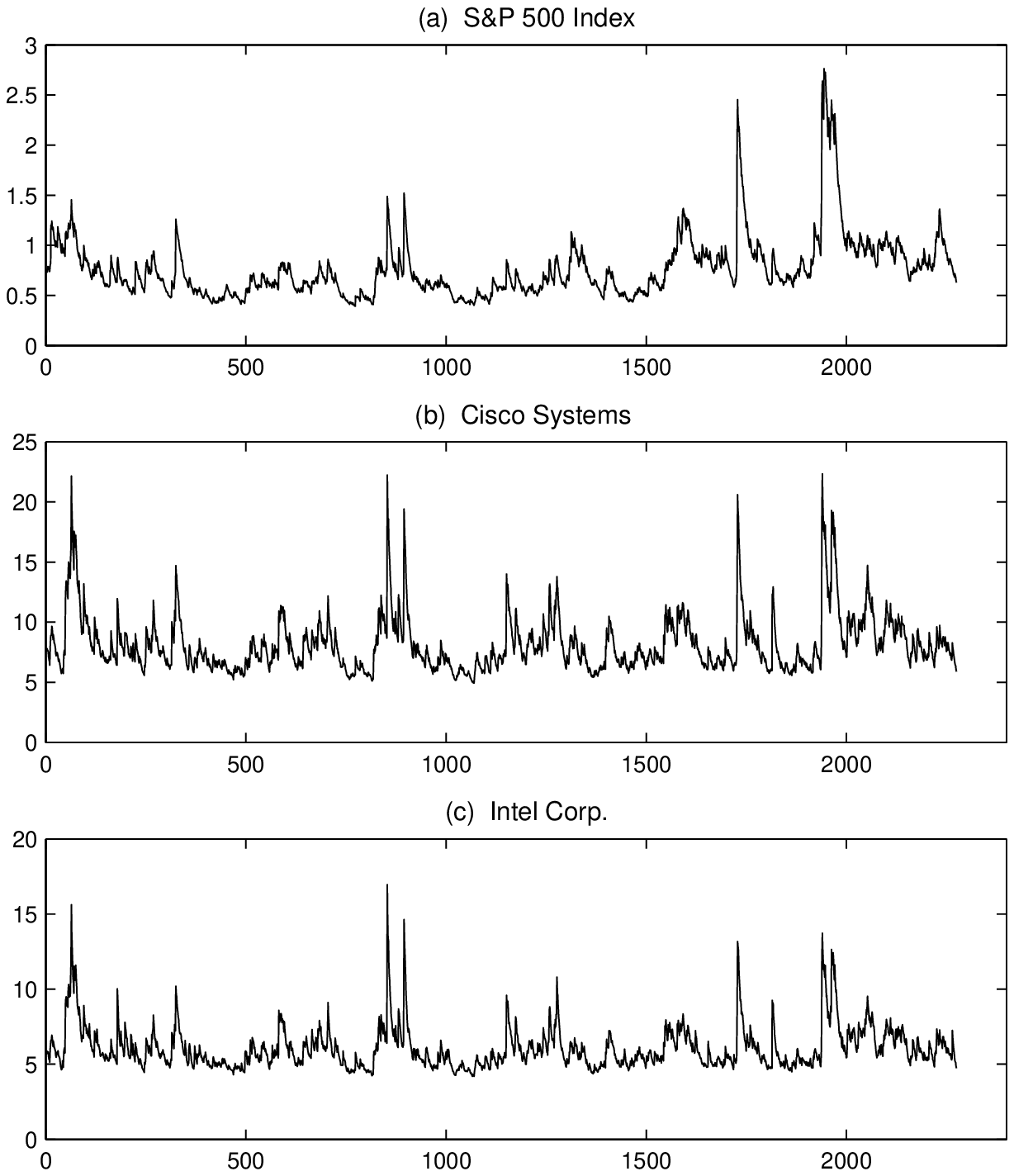}}
\begin{singlespace}

\caption[Fig 5] {\sl Fitted volatility processes based on
Orthogonal-GARCH(1,1) model for daily log returns of (a)  $S\&P$
500 index, (b) Cisco Systems stock  and (c) Intel Corporation
stock. }

\end{singlespace}
\end{figure}

\newpage

\begin{figure}
\centerline{\psfig{figure=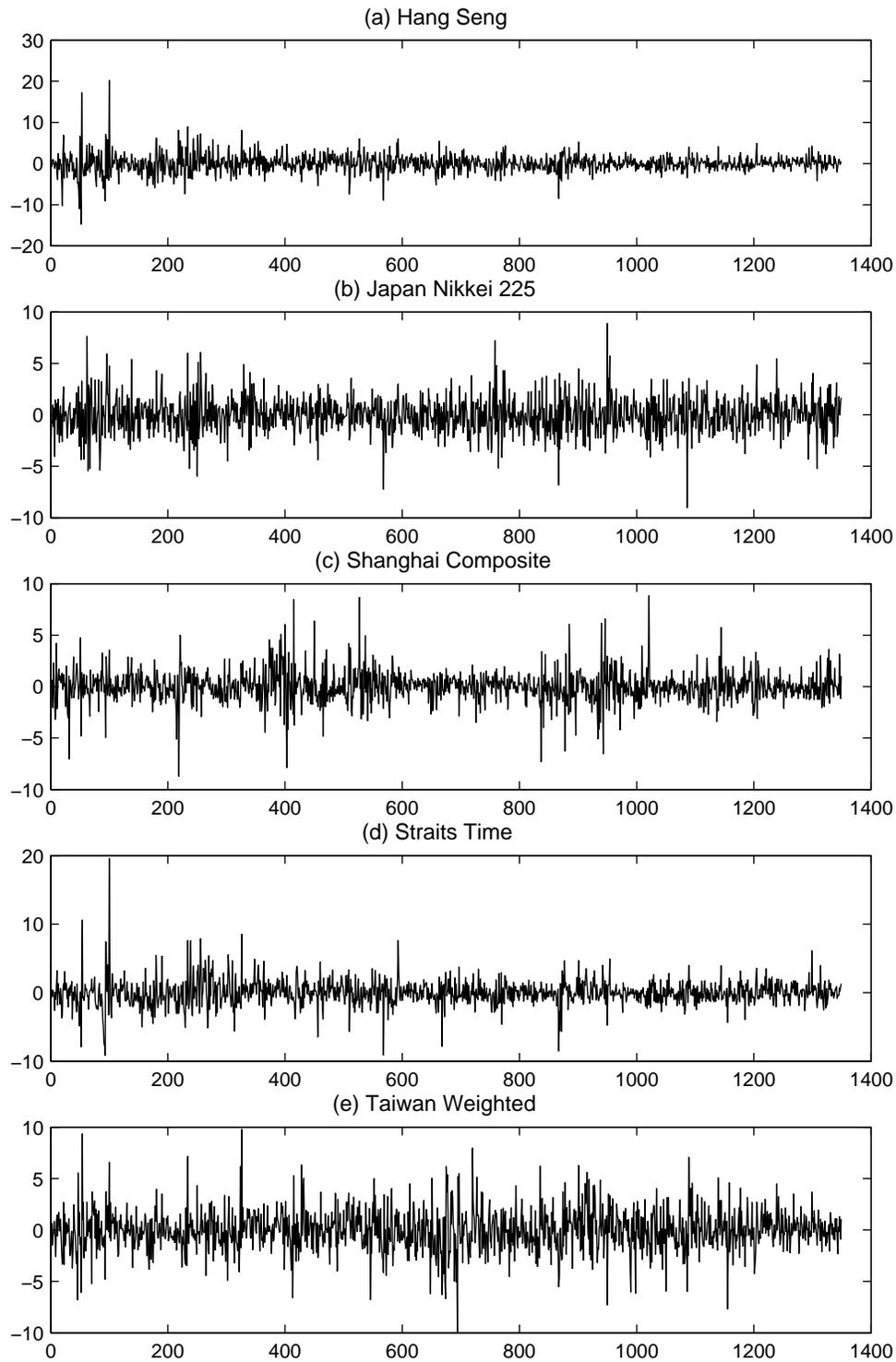}}
\begin{singlespace}

\caption[Fig 6] {\sl Plots of dividend adjusted daily log returns of
(a)  Hang Seng index in Hong Kong, (b) Japan Nikkei 225 index, (c)
Shanghai Composite index in China, (d) Singapore Straits Time index,
 and (e) Taiwan Weighted index. Time span is from August 1, 1997 to December 30, 2003 with 1349 observations.}

\end{singlespace}
\end{figure}

\newpage

\begin{figure}
\centerline{\psfig{figure=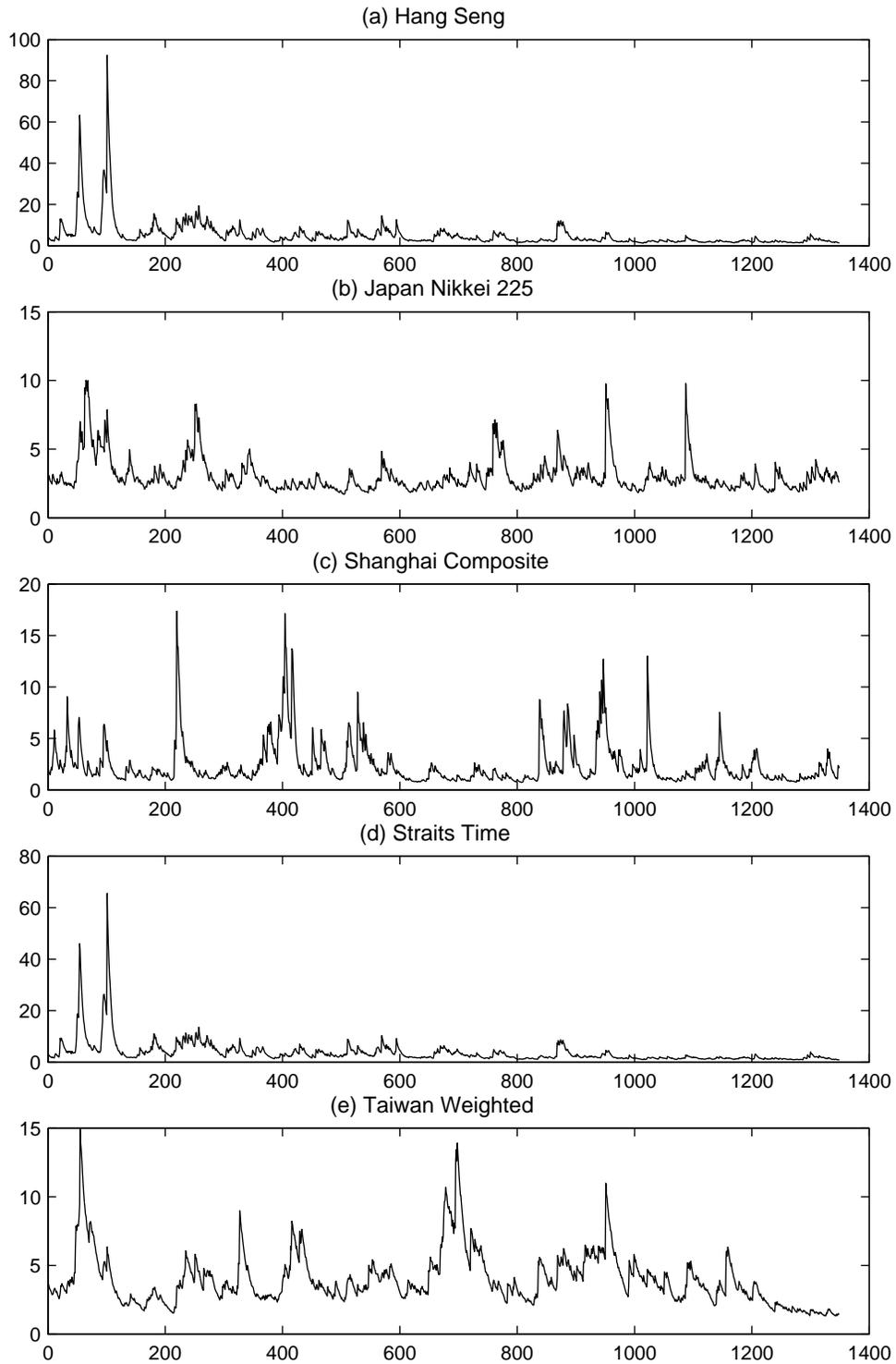}}
\begin{singlespace}

\caption[Fig 7] {\sl Fitted volatility processes based on CUC-Extended GARCH(1,1) model for daily log returns of
(a)  Hang Seng index in Hong Kong, (b) Japan Nikkei 225 index, (c)
Shanghai Composite index in China, (d) Singapore Straits Time index,
 and (e) Taiwan Weighted index. }

\end{singlespace}
\end{figure}

\newpage

\begin{figure}
\centerline{\psfig{figure=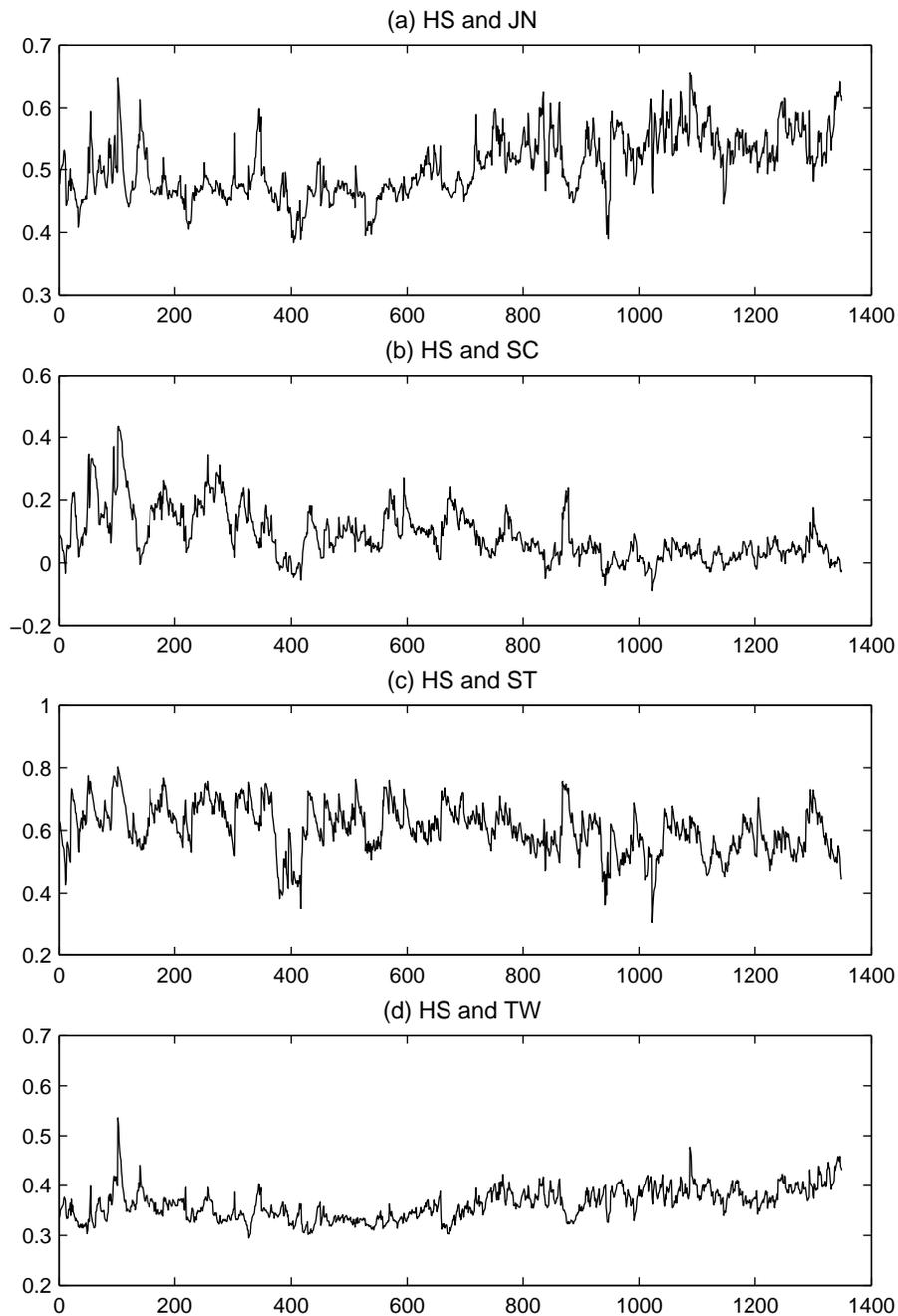}}
\begin{singlespace}

\caption[Fig 7] {\sl Fitted conditional correlations between daily
log-returns of Hang Seng index (HS) and (a) Japan Nikkei 225 index (JN),
(b) Shanghai
Composite index in China (SC), (c) Singapore Straits Time index (ST),
 (d) Taiwan Weighted index (TW).}

\end{singlespace}
\end{figure}

\end{document}

%% file: QYalias.tex
\@addtoreset{equation}{section}
\renewcommand{\theequation}{\thesection.\arabic{equation}}
\renewcommand{\hat}{\widehat}
\def\singlespace{\def\baselinestretch{1}\@normalsize}

\def\wh{\widehat}
\def\wt{\widetilde}
\def\askip{\vspace{0.1in}}

\newcommand{\cov}{{\rm Cov}}
\newcommand{\diag}{{\rm diag}}

\newcommand{\sgn}{\mbox{sgn}}

\newcommand{\tr}{\mbox{tr}}
\newcommand{\var}{\mbox{Var}}
\def\ga{\gamma}

\def\la{\lambda}

\newcommand{\ve}{{\varepsilon}}

\newcommand{\bA}{{\mathbf A}}
\newcommand{\bB}{{\mathbf B}}

\newcommand{\bE}{{\mathbf E}}

\newcommand{\bI}{{\mathbf I}}

\newcommand{\bR}{{\mathbf R}}
\newcommand{\bS}{{\mathbf S}}
\newcommand{\bU}{{\mathbf U}}
\newcommand{\bV}{{\mathbf V}}

\newcommand{\bX}{{\mathbf X}}
\newcommand{\bY}{{\mathbf Y}}
\newcommand{\bZ}{{\mathbf Z}}
\newcommand{\ba}{{\mathbf a}}
\newcommand{\bb}{{\mathbf b}}

\newcommand{\bu}{{\mathbf u}}
\newcommand{\bv}{{\mathbf v}}

\newcommand{\bvarphi}  {\boldsymbol{\varphi}}

\newcommand{\bSigma}{\boldsymbol{\Sigma}}
\newcommand{\bDelta}{\boldsymbol{\Delta}}
\newcommand{\bgamma}{\boldsymbol{\gamma}}

\newcommand{\bve}{\mbox{\boldmath$\varepsilon$}}

\newcommand{\btheta} {\boldsymbol{\theta}}

\newcommand{\bLambda} {\boldsymbol{\Lambda}}
\newcommand{\bC}{{\mathbf C}}
\newcommand{\bD}{{\mathbf D}}

\newcommand{\scon}{\stackrel{a.s.}{\longrightarrow}}

\newcommand{\calF}{{\mathcal F}}
\newcommand{\calH}{{\mathcal H}}

\newcommand{\calB}{{\mathcal B}}

\newcommand{\RR}{\EuScript R}